\documentclass[11pt,a4paper]{article}
\usepackage[left=25mm,right=25mm,top=19mm,bottom=19mm,footskip=10mm]{geometry}
\usepackage{amsmath,amsthm,amscd,amssymb}
\usepackage{txfonts}
\usepackage{bm}
\usepackage{mathtools}
\usepackage{mathrsfs}
\usepackage{indentfirst}
\usepackage{titlesec}
\usepackage{microtype}
\usepackage{hyperref}
\hypersetup{colorlinks=true,linkcolor=blue,citecolor=blue,urlcolor=blue}
\allowdisplaybreaks[4]
\titleformat{\section}{\normalfont\normalsize\bfseries}{\thesection.}{0.5em}{}
\titlespacing*{\section}{0pt}{2.5ex plus 1ex minus 0.2ex}{1.1ex plus 0.2ex}
\titleformat{\subsection}{\normalfont\normalsize\itshape}{\thesubsection.}{0.5em}{}
\titlespacing*{\subsection}{0pt}{1.8ex plus 0.8ex minus 0.2ex}{0.6ex plus 0.2ex}

\newtheorem{thm}{Theorem}[section]
\newtheorem{cor}[thm]{Corollary}
\newtheorem{lem}[thm]{Lemma}
\newtheorem{prop}[thm]{Proposition}
\theoremstyle{definition}
\newtheorem{defn}[thm]{Definition}
\theoremstyle{remark}
\newtheorem{rem}[thm]{Remark}
\newtheorem{exm}[thm]{Example}
\numberwithin{equation}{section}
\newcommand{\pf}{\begin{proof}}
\newcommand{\zb}{\end{proof}}

\newcommand{\ran}{\mathop{\rm ran}\nolimits}

\begin{document}
\thispagestyle{empty}
\begin{center}
{\Large Almost invariant subspaces of the Ces\`aro operator and nearly invariant semigroup models\par}
\vspace{1.45em}
{\normalsize Yuxia Liang\par}

\end{center}
 
\noindent\rule{\linewidth}{0.4pt}\par
\vspace{0.5em}
\noindent\textbf{Abstract}\par\smallskip
\noindent

We investigate the finite Blaschke extensions
$K_{B\Theta_\delta}$ of Ces\`aro-invariant model
spaces in $H^2(\mathbb D)$, where $B$ is a finite
Blaschke product and $\Theta_\delta$ is an atomic
singular inner function. We prove that the
Ces\`aro invariance defect of $K_{B\Theta_\delta}$
equals the number of distinct zeros of $B$,
independently of their multiplicities.
For nonconstant $B$, the largest Ces\`aro-invariant
subspace contained in $K_{B\Theta_\delta}$ is
$K_{\Theta_\delta}$, while the smallest closed
Ces\`aro-invariant subspace containing it is $H^2$.
Moreover, $K_{B\Theta_\delta}$ is nearly invariant
for the associated composition semigroup if and
only if $B$ is a unimodular multiple of a power
of $z$. These results are obtained through shift
and weighted-translation models, which also yield
nearly invariant extensions, a model-space rigidity
theorem, and a general exact defect formula
governed by the successor relation
$\lambda\mapsto\lambda+1$.

\par
\vspace{0.55em}
\noindent\rule{\linewidth}{0.4pt}\par

\section{Introduction}\label{s1}

The Ces\`aro operator $\mathcal C$ on $H^2:=H^2(\mathbb D)$ is given by
\[
(\mathcal Cf)(z)=\frac1z\int_0^z\frac{f(\zeta)}{1-\zeta}\,d\zeta,
\]
where the value at $z=0$ is defined by continuity. Its classical operator
theory goes back to Brown, Halmos and Shields \cite{BrownHalmosShields},
Kriete and Trutt \cite{KrieteTrutt,KrieteTrutt1975}, and Cowen
\cite{Cowen1984}. The connection with composition semigroups, developed by
Siskakis \cite{Siskakis1987} and subsequently by
Gallardo-Guti\'errez and Partington \cite{GallardoPartingtonCesaro} and Gallardo-Guti\'errez, Partington and Ross
\cite{GallardoPartingtonRoss},
provides useful models for its invariant-subspace problem.

For an inner function $u$, let $K_u:=H^2\ominus uH^2$ and write the atomic singular inner function
\begin{equation*} 
\Theta_\delta(z):=\exp\!\left(-\delta\frac{1+z}{1-z}\right),
\qquad \delta>0.
\end{equation*}
Recent work of Gallardo-Guti\'errez, Partington and Ross
\cite{GallardoPartingtonRoss} shows that the nonzero model spaces invariant
under $\mathcal C$ are precisely $K_{\Theta_\delta}$ and describes their
contained invariant subspaces. In particular, their Theorem~6.3 shows that
every proper closed $\mathcal C$-invariant subspace containing
$K_{\Theta_\delta}$ is $K_{\Theta_\beta}$ for some $\beta\ge\delta$.
Despite these advances, a complete description of the closed invariant subspaces of the Ces\`aro operator remains unavailable. This motivates the study of almost invariant subspaces, which provide a natural finite-dimensional relaxation of ordinary invariance.

Throughout the paper, $\bigvee$ denotes closed linear span. For a
bounded operator $T$ and a subset $A$, put
$[A]_T:=\bigvee\{T^nf:f\in A,\ n\ge0\}$; for a closed subspace
$M\subseteq H^2$, define
\[
\mathcal I_{\mathcal C}(M)
:=\{f\in M:\mathcal C^nf\in M\text{ for all }n\ge0\}.
\]
Thus $\mathcal I_{\mathcal C}(M)$ is the largest closed
$\mathcal C$-invariant subspace contained in $M$, whereas
$[M]_{\mathcal C}$ is the smallest closed $\mathcal C$-invariant
subspace containing $M$.
For a bounded operator $T$ and a closed subspace $M$, define its
\emph{invariance defect} by
$$\operatorname{def}_T(M):=\operatorname{rank}(P_{M^\perp}T|_M).$$
The subspace $M$ is $T$-invariant if and only if
$\operatorname{def}_T(M)=0$, and is \emph{almost invariant}
if $\operatorname{def}_T(M)<\infty$, or equivalently,
if $TM\subseteq M+F$ for some finite-dimensional subspace $F$.
Recall that a \emph{half-space} is a closed subspace of
infinite dimension and infinite codimension.
Since every closed subspace that is not a half-space
is almost invariant under any bounded operator,
we restrict our attention to almost invariant half-spaces
throughout the remainder of this paper.
Almost-invariant subspaces have been studied in connection with the
invariant-subspace problem, finite-rank perturbations, and operator
algebras; see Marcoux, Popov and Radjavi \cite{MarcouxPopovRadjavi},
Popov and Tcaciuc \cite{PopovTcaciuc}, and Sirotkin and Wallis
\cite{SirotkinWallis}. For Hardy-space shifts, Chalendar,
Gallardo-Guti\'errez and Partington obtained Beurling-type descriptions
\cite{ChalendarGallardoPartington}; more recently, Gu and collaborators
related almost-invariant shift subspaces to Toeplitz and Hankel operators
\cite{GuHwangKimLeePark}. These results provide context for the exact
Ces\`aro defects considered in this paper.

For a finite Blaschke product $B$ of degree $N$, the decomposition
\begin{equation}\label{e2}
K_{B\Theta_\delta}
=K_{\Theta_\delta}\oplus\Theta_\delta K_B,
\qquad \dim K_B=N,
\end{equation}
 shows that the half-space $K_{B\Theta_\delta}$ is
a finite-dimensional extension of the
$\mathcal C$-invariant model space $K_{\Theta_\delta}$.
Consequently,
$\operatorname{def}_{\mathcal C}K_{B\Theta_\delta}\le N$.
This observation gives no information about the exact
defect or its dependence on the zeros of $B$.
We therefore ask for the following main questions.  
\vspace{1mm}

{\bf Main\;\;Questions}: What is the exact defect of $K_{B\Theta_\delta}$, how does it depend on the zeros and multiplicities of $B$, and how is the extension related to genuine invariant subspaces?

Another related question is which of these almost invariant
extensions are nearly invariant for the associated
composition semigroup. To investigate these questions, we use the connection between the Ces\`aro operator and the associated affine composition semigroup induced by
\begin{equation}\label{e3}
\varphi_t(z):=e^{-t}z+1-e^{-t},
\qquad C_tf:=f\circ\varphi_t,\qquad t\ge0.
\end{equation}
A closed subspace is invariant under $\mathcal C$ if and only if its
orthogonal complement is invariant under every $C_t$
\cite[Theorem~2.1]{GallardoPartingtonCesaro}. Moreover,
\begin{equation}\label{e4}
\mathcal C=\int_0^\infty e^{-t}C_t^*\,dt,
\end{equation}
where norm convergence is established below in Lemma~\ref{r4}.

To study the main questions and corresponding near invariance
for the composition semigroup, we need the
weighted-translation model on the weighted space
\begin{equation}\label{e6}
L_w^2:=L^2(\mathbb R,w(y)\,dy),
\qquad w(y):=e^{-2(e^y-1)}.
\end{equation}
For $a\in\mathbb R$, put
\begin{equation}\label{e7}
\mathcal N_a:=L^2((-\infty,a),w(y)\,dy),
\qquad \delta:=e^a.
\end{equation}
The unitary map $\mathcal F$ in \eqref{e14} sends $K_{\Theta_\delta}$
onto $\mathcal N_a$ and intertwines
$\{C_t\}_{t\ge0}$ with $\{e^{t/2}S_t\}_{t\ge0}$,
where $S_t$ is right translation on $L_w^2$,
defined in \eqref{STT}.
Thus $\mathcal N_a$ provides the natural base
space for studying the finite-dimensional
extensions $K_{B\Theta_\delta}$.  Alongside almost invariance for $\mathcal C$, we define the following
preimage-based notion for semigroups.

\begin{defn}\label{r1}
Let $\{T_t\}_{t\ge0}$ be a $C_0$-semigroup on a Hilbert space
$\mathcal H$. A closed subspace $M\subseteq\mathcal H$ is
\emph{nearly invariant for $\{T_t\}_{t\ge0}$} if
\begin{equation}\label{e5}
T_tf\in M\text{ for some }t>0
\quad\Longrightarrow\quad f\in M.
\end{equation}
Equivalently, $T_t^{-1}(M)\subseteq M$ for every $t>0$, where
$T_t^{-1}(M):=\{f\in\mathcal H:T_tf\in M\}$ is a set-theoretic
preimage.
\end{defn}

Following \cite[Definition~1.4]{LP22}, this is also called near
$\{T_t^*\}_{t\ge0}$-invariance. When $T_t$ is an isometry,
\eqref{e5} is equivalent to
$T_t^*(M\cap\operatorname{ran}T_t)\subseteq M$ for each $t>0$. The related theory of nearly backward-shift-invariant subspaces and
Toeplitz kernels provides further motivation in
\cite{Hayashi,Hitt,SarasonNear,SarasonToeplitz,LP21}.
However, for arbitrary semigroups it is different from the invariance
under the adjoints $T_t^*$. Indeed,
$T_t^{-1}(M)=(T_t^*M^\perp)^\perp$ and hence \eqref{e5} is
equivalent to $M^\perp\subseteq\overline{T_t^*M^\perp}$ for every
$t>0$. In particular, if $T_tf\in M$, the same condition applied to
$T_{t-s}(T_sf)=T_tf$ gives $T_sf\in M$ whenever $0\le s\le t$.

To address these questions, we use two complementary
semigroup models. The isometric-shift model $P_t$ on
$L^2(0,\infty)$, defined in \eqref{Pt}, yields
minimal nearly invariant subspaces through backward
orbits; the corresponding rational-orbit constructions
are related to earlier work
\cite{LP22,LiangPartingtonCyclic}.
In the weighted-translation model, the domains of  backward
translations $\mathcal{B}_t$ are   taken into account and the  criterion for
near invariance is presented in Theorem~\ref{r7}. Within the exponential--polynomial
extensions arising from this model,
\emph{Theorem~\ref{r17} characterizes exactly which spaces
are Hardy model spaces, showing that the exponential
rates and their multiplicities are rigidly prescribed.}

 Our main defect calculations rely on the integral
representation of $\mathcal C$ and its differential
inverse. These turn the defect problem for
exponential--polynomial tails into a finite-dimensional
calculation governed by the successor relation
$\lambda\mapsto\lambda+1$ (Theorem~\ref{r24}). \emph{
For a finite Blaschke product $B$,
Theorem~\ref{r27} gives
\[
\operatorname{def}_{\mathcal C}K_{B\Theta_\delta}
=\#\{z\in\mathbb D:B(z)=0\},
\]
where zeros are counted without multiplicity.}
In particular, $K_{z^N\Theta_\delta}$ has defect one,
although it is an $N$-dimensional extension of
$K_{\Theta_\delta}$.
Using the known classification of closed
$\mathcal C$-invariant subspaces containing
$K_{\Theta_\delta}$, Proposition~\ref{r22}
identifies $K_{\Theta_\delta}$ as the largest
invariant subspace contained in $K_{B\Theta_\delta}$
and $H^2$ as the smallest closed invariant subspace
containing it, for nonconstant $B$.
Finally, Theorem~\ref{r30} shows that
$K_{B\Theta_\delta}$ is nearly invariant for
$\{C_t\}_{t\ge0}$ precisely when $B$ is a power
of $z$, up to a unimodular constant.
Thus almost invariance for $\mathcal C$ and
near invariance for its composition semigroup
are distinct conditions.

So the paper is organized as follows. Section~2 establishes the semigroup and integral
models. Section~3 constructs nearly invariant
extensions and proves the model-space rigidity
theorem. Section~4 derives the exact Ces\`aro
defects, determines invariant hulls, and
characterizes near invariance for finite
Blaschke extensions.

\section{Semigroup and integral models}\label{s2}
This section establishes the unitary model for the affine semiflow
\eqref{e3}. For background on composition semigroups, see
\cite{AvicouChalendarPartington,AvicouChalendarPartington2,BerksonPorta,BracciContrerasDiaz,CowenMacCluer,EngelNagel,GallardoYakubovich,SiskakisReview};
for translation models and invariant subspaces, see
\cite{Domar1981,Domar1983,GallardoPartington2009,GallardoPartingtonRodriguezMathZ,Lax}.
Write
$ \mathbb C_+:=\{s\in\mathbb C:\Re s>0\}.$ The Hardy-space norms are fixed by
\[
 \|\sum_{n\geq0}a_nz^n \|_{H^2(\mathbb D)}^2
=
\sum_{n\geq0}|a_n|^2
\]
and
\[
\|G\|_{H^2(\mathbb C_+)}^2
=
\frac1{2\pi}
\sup_{\sigma>0}
\int_{\mathbb R}|G(\sigma+i\xi)|^2\,d\xi.
\]
 Let
\begin{equation}\label{e8}
b(s):=\frac{s-1}{s+1},
\qquad
q(z):=\frac{1+z}{1-z},
\end{equation}
so that $b:\mathbb C_+\to\mathbb D$ and
$q:\mathbb D\to\mathbb C_+$ are inverse conformal maps. For every analytic
function $u$ on $\mathbb D$, the superscript $+$ denotes its
right-half-plane counterpart $u^+:=u\circ b$. In particular,
$\Theta_\delta^+(s)=e^{-\delta s}$. Define
\begin{equation*}
(Uf)(s)
:=
\frac{\sqrt2}{s+1}f(b(s)),
\qquad
(U^{-1}G)(z)
=
\frac{\sqrt2}{1-z}G(q(z)).
\end{equation*}
With the norms above, $U:H^2(\mathbb D)\to H^2(\mathbb C_+)$ is
unitary. Moreover, if $\Theta$ is inner on $\mathbb C_+$, then
\begin{equation}\label{e9}
U^{-1}K_\Theta
=
K_{\Theta\circ q}.
\end{equation}

A direct calculation gives
 $ q(\varphi_t(z)) =e^tq(z)+(e^t-1).$ 
Consequently,
\begin{equation}\label{e10}
UC_tU^{-1}
=
e^t C_{\phi_t},
\qquad
\phi_t(s):=e^ts+e^t-1.
\end{equation}
Let $\mathcal L: L^2(0, \infty)\to H^2(\mathbb C_+)$ be the Laplace
transform
\[
(\mathcal Lg)(s)
:=
\int_0^\infty g(x)e^{-sx}\,dx,
\qquad g\in L^2(0, \infty).
\]
By the Paley--Wiener theorem and the chosen normalization,
$\mathcal L $ is unitary.
Let $  T:L^2(\mathbb R)\to L^2(0, \infty)$ be
\[
(  Th)(x):=x^{-1/2}h(\log x),
\qquad x>0,
\]
with inverse
\begin{equation}\label{e11}
(  T^{-1}g)(y)=e^{y/2}g(e^y).
\end{equation}
The standard computation for \eqref{e10}
gives
\[
\mathcal L^{-1}C_{\phi_t}\mathcal L=V_t,
\qquad
(V_tg)(x)
=
e^{-t}e^{-(1-e^{-t})x}g(e^{-t}x),
\]
and
\[
  T^{-1}V_t  T
=
e^{-t/2}\sigma_t,
\qquad
(\sigma_th)(y)
=
e^{-(1-e^{-t})e^y}h(y-t).
\]
On the weighted space $L_w^2$ defined in
\eqref{e6},    the unitary operator $W: L^2(\mathbb R)\rightarrow L_w^2$ is defined by
\begin{equation}\label{e12}
(Wh)(y):=\frac{h(y)}{\sqrt{w(y)}}.
\end{equation}
As in \cite[Proposition~2.4]{GallardoPartingtonCesaro},
$W\sigma_tW^{-1}=S_t$, where
\begin{align}\label{STT}
(S_tg)(y):=g(y-t),\qquad t\geq0. 
\end{align}
Thus $S_t$ is right translation acting on the weighted space $L_w^2$.
The family $\{S_t\}_{t\geq0}$ is a $C_0$-semigroup of contractions.
Indeed, since
\[
w(y+t)
=e^{-2(e^{y+t}-1)}
\leq e^{-2(e^y-1)}
=w(y),
\qquad t\geq0,
\]
a change of variables gives
\begin{align}\label{e13}
\|S_tg\|_{L_w^2}^2
=
\int_{\mathbb R}|g(y)|^2w(y+t)\,dy
\leq
\int_{\mathbb R}|g(y)|^2w(y)\,dy
=
\|g\|_{L_w^2}^2.
\end{align}
In fact, $\|S_t\|=1$. To see this,
\[
\|S_t\|^2
=
\operatorname*{ess\,sup}_{y\in\mathbb R}
\frac{w(y+t)}{w(y)}.
\]
Since \(w(y+t)\le w(y)\), the right-hand side is at most \(1\), whereas
$
\frac{w(y+t)}{w(y)}\longrightarrow1,$ as $
  y\to-\infty, $ 
shows that its essential supremum is at least $1$. Thus
$\|S_t\|=1$. Strong continuity follows first on
$C_c^\infty(\mathbb R)$ from ordinary translation continuity and
local equivalence of $w(y)\,dy$ with Lebesgue measure, and then
follows on all of $L_w^2$ from density and
\eqref{e13}. 

Finally, the unitary operator
\begin{equation}\label{e14}
\mathcal F:=WT^{-1}\mathcal L^{-1}U:
H^2(\mathbb D)\longrightarrow L_w^2,
\end{equation}
therefore gives the intertwining identity
\begin{equation}\label{e15}
\mathcal F C_t\mathcal F^{-1}
=T_t:=e^{t/2}S_t,
\qquad t\ge0.
\end{equation} 
In particular, $\|T_t\|=e^{t/2}$, whereas $\|S_t\|=1$.
Since $e^{t/2}\ne0$,   $\{T_t\}_{t\geq 0}$ and $\{S_t\}_{t\geq 0}$
have exactly the same nearly invariant subspaces. Under $\mathcal F$,
near invariance for $\{C_t\}_{t\geq 0}$ is therefore equivalent to that for
$\{S_t\}_{t\geq 0}$.
\begin{lem}\label{r2}
Let $\mathcal F$ be defined by \eqref{e14}. A closed subspace
$\mathcal M\subseteq H^2(\mathbb D)$ is nearly invariant for
$\{C_t\}_{t\geq0}$ if and only if $\mathcal F\mathcal M$
is nearly invariant for $\{S_t\}_{t\geq0}$ on $L_w^2$.
Equivalently, $\mathcal F\mathcal M$ is nearly invariant for
$\{T_t\}_{t\geq0}$.
\end{lem}

\begin{proof}
Suppose that $\mathcal M$ is nearly invariant for $\{C_t\}_{t\ge0}$,
and let $S_tg\in\mathcal F\mathcal M$ for some $t>0$.
Because $T_tg=e^{t/2}S_tg$, equation \eqref{e15} gives
$C_t\mathcal F^{-1}g\in\mathcal M$. Hence
$\mathcal F^{-1}g\in\mathcal M$ and $g\in\mathcal F\mathcal M$.
Conversely, if $C_tf\in\mathcal M$, then
$S_t\mathcal Ff=e^{-t/2}\mathcal FC_tf\in\mathcal F\mathcal M$.
Near invariance for $\{S_t\}_{t\ge0}$ implies $\mathcal Ff\in\mathcal F\mathcal M$,
and hence $f\in\mathcal M$. The equivalence with $\{T_t\}_{t\ge0}$
follows from $T_t=e^{t/2}S_t$.
\end{proof}

We next record the adjoint of $S_t$ on $L_w^2$, whose formula involves the weight.

\begin{lem}\label{r3}
For every $t\geq0$, the adjoint of the right-translation operator
$S_t$ on $L_w^2$ is
\begin{equation}\label{e16}
(S_t^*g)(y)
=
\frac{w(y+t)}{w(y)}g(y+t)
=
e^{-2e^y(e^t-1)}g(y+t), \qquad y\in \mathbb{R}.
\end{equation}
Moreover, $\{S_t^*\}_{t\geq0}$ is a $C_0$-semigroup of
contractions.
\end{lem}

\begin{proof}
For $f,g\in L_w^2$,
\[
\begin{aligned}
\langle S_tf,g\rangle_{L_w^2}
&=
\int_{\mathbb R}f(y-t)\overline{g(y)}w(y)\,dy\\
&=
\int_{\mathbb R}
f(x)\overline{g(x+t)}w(x+t)\,dx,
\end{aligned}
\]
which proves \eqref{e16}. The contraction
property follows from $\|S_t^*\|=\|S_t\|\leq1$. Since
$\{S_t\}_{t\geq 0}$ is a $C_0$-semigroup on  $L_w^2$, its adjoint
semigroup is strongly continuous; alternatively,
\[
\|S_t^*f-f\|^2
\leq
2\|f\|^2-2\Re\langle S_t^*f,f\rangle
\longrightarrow0, 
\qquad t\rightarrow 0. \qedhere
\]
\end{proof}

The invariant subspaces of translation semigroups
on $L^2(\mathbb R)$ are described by classical
Beurling--Lax theory; see \cite{Lax}. In contrast,
the invariant-subspace structure of the weighted
right-translation semigroup $\{S_t\}_{t\geq0}$
on $L_w^2$ is more complicated; see
\cite{Domar1981,Domar1983,GallardoPartingtonCesaro}.
We therefore begin with  \emph{ the support subspaces
$\mathcal N_a$ defined in \eqref{e7}}, which provide
a simple setting for studying both invariance
and near invariance.  For the near invariance,  if $S_tg\in\mathcal N_a$, then
$g$ vanishes a.e. on $[a-t,\infty)$ and hence belongs to
$\mathcal N_a$. 
Since $W^{-1}$ preserves supports and
$T$ transforms the interval $(-\infty,a)$
onto $(0,e^a)$, we have
\[
W^{-1}\mathcal N_a=L^2(-\infty,a),
\qquad
TW^{-1}\mathcal N_a
=L^2(0,e^a)=L^2(0,\delta).
\]
By the Paley--Wiener theorem and \eqref{e9},
\begin{align}\label{e17}
\mathcal L TW^{-1}\mathcal N_a
 =\mathcal L L^2(0,\delta)
=K_{e^{-\delta s}}, \qquad
\mathcal F^{-1}\mathcal N_a
 =U^{-1}K_{e^{-\delta s}}
=K_{\Theta_\delta},
\end{align}
where $\Theta_\delta$ is given by \eqref{e2}. The adjoint orbit can also be identified exactly 
\begin{equation}\label{e18}
S_t^*\mathcal N_a=\mathcal N_{a-t},
\qquad
C_t^*K_{\Theta_\delta}
= K_{\Theta_{\delta e^{-t}}},
\qquad t\geq0.
\end{equation}
Indeed, the inclusion
$S_t^*\mathcal N_a\subseteq\mathcal N_{a-t}$ follows from
\eqref{e16}. Conversely, if
$h\in\mathcal N_{a-t}$, define
\[
g(x)
=
h(x-t)\frac{w(x-t)}{w(x)},
\qquad x<a,
\]
and $g(x)=0$ for $x\geq a$. On $(-\infty,a)$, since 
\[
\frac{w(x-t)}{w(x)}
=
\exp\!\bigl(2e^x(1-e^{-t})\bigr)
\le
\exp\!\bigl(2e^a(1-e^{-t})\bigr),
\qquad x<a,
\] $g\in\mathcal N_a$ and $S_t^*g=h$. The disk identity in
\eqref{e18} follows from
\eqref{e15}, \eqref{e17}, and
$\delta e^{-t}=e^{a-t}$. Then \eqref{e18} 
implies  $K_{\Theta_\delta}$ is
$\{C_t^*\}_{t\geq0}$-invariant. Its near invariance for $\{C_t\}_{t\geq0}$ follows from
Lemma~\ref{r2}.

\subsection{An isometric shift and the Ces\`aro integral model}\label{s3}
With the unitary maps $U$ and $\mathcal L$ from
Section~\ref{s2}, define the unitary map
\begin{align}
\mathcal J:=\mathcal L^{-1}U:H^2(\mathbb D)\longrightarrow L^2(0,\infty).\label{JJ}
\end{align}
For an inner function $u$, \eqref{e9} gives $UK_u=K_{u^+}$.

To construct minimal nearly invariant subspaces,  we use the
classical isometric right-shift semigroup $\{P_t\}_{t\ge0}$ on $L^2(0,\infty)$, which is defined by
\begin{equation}\label{Pt}
(P_tf)(x):=\chi_{(t,\infty)}(x)f(x-t),\qquad t\ge0,
\end{equation}
for $f\in L^2(0,\infty)$. Each $P_t$ defined in \eqref{Pt} is an isometry and
\begin{align}\mathcal LP_t\mathcal L^{-1}=M_{e^{-ts}}.\label{mus} \end{align}
Related properties have been studied in \cite{LP22,LiangPartingtonCyclic}.

To represent the adjoint composition semigroup, we use the same
Cayley--Laplace transform $\mathcal J$ in \eqref{JJ}. The calculation in
\cite[Section 2]{GallardoPartingtonCesaro} gives
\begin{equation}\label{e19}
Q_t:=\mathcal J C_t^*\mathcal J^{-1},\qquad
(Q_tf)(x)=e^t e^{-(e^t-1)x}f(e^t x),\qquad t\ge0.
\end{equation}
For completeness, the identity \eqref{e10} gives
$UC_tU^{-1}=e^t C_{\phi_t}$ with
$\phi_t$ as defined  in \eqref{e10}, and Laplace inversion gives the corresponding
operator on $L^2(0,\infty)$; taking its adjoint yields \eqref{e19}.
Thus $Q_t$ represents $C_t^*$ on $L^2(0,\infty)$.

The proof of \cite[Theorem~2.1]{GallardoPartingtonCesaro} gives
$\mathcal C^*=(I-A)^{-1}$ for the composition-semigroup generator $A$.
The next lemma records the resulting integral model with our normalization.

\begin{lem}\label{r4}
The integral $\mathcal C=\int_0^\infty e^{-t}C_t^*\,dt$
converges in operator norm. Under the unitary map $\mathcal J$ in \eqref{JJ},
the Ces\`aro operator is represented by
\begin{equation}\label{e20}
V:=\mathcal J{\mathcal C}\mathcal J^{-1},\qquad
(Vf)(x)=e^x\int_x^\infty e^{-u}\frac{f(u)}{u}\,\mathrm d u,
\quad x>0.
\end{equation}
If $h=Vf$, then $h$ is locally absolutely continuous on $(0,\infty)$ and
\begin{equation}\label{e21}
f(x)=x(h(x)-\frac{\mathrm d h}{\mathrm d x}(x)) 
\qquad\text{for a.e. }x>0.
\end{equation}
\end{lem}
\begin{proof}
By \eqref{e15}, we have $\|C_t\|=e^{t/2}$ for
$t\ge0$. Hence, for $0\le R<S$,
\[
\left\|\int_R^S e^{-t}C_t\,dt\right\|
\leq\int_R^S e^{-t}\|C_t\|\,dt
\leq 2e^{-R/2}.
\]
Thus the strong integrals over finite intervals
converge in operator norm as $S\to\infty$.
Since the growth bound of $\{C_t\}_{t\ge0}$
is at most $\frac{1}{2}$, the semigroup resolvent
formula applies at $\lambda=1$ and yields
\[
(I-A)^{-1}=\int_0^\infty e^{-t}C_t\,dt.
\]
Combining this with
$\mathcal C^*=(I-A)^{-1}$ gives the desired
representation of $\mathcal C^*$. Taking
adjoints is justified by operator-norm
convergence and gives \eqref{e4}. Moreover, 
$\|C_t^*\|=\|C_t\|$.

We next derive the integral representation
of $V=\mathcal J\mathcal C\mathcal J^{-1}$.
By \eqref{e19},
\[
(\mathcal JC_t^*\mathcal J^{-1}f)(x)
=e^t e^{-(e^t-1)x}f(e^t x).
\]
Consequently,
\[
(Vf)(x)=\int_0^\infty
e^{-t}e^t e^{-(e^t-1)x}f(e^t x)\,dt=e^x\int_x^\infty
e^{-u}\frac{f(u)}{u}\,du.
\]

To justify the integral representation
for arbitrary $f\in L^2(0,\infty)$,
write $T_t=\mathcal JC_t^*\mathcal J^{-1}$.
For every $0<a<b<\infty$, the
Cauchy--Schwarz inequality gives
\[
\int_a^b\int_0^\infty
e^{-t}|(T_tf)(x)|\,dt\,dx\leq (b-a)^{\frac{1}{2}}
\int_0^\infty e^{-t}\|T_tf\|_2\,dt\leq 2(b-a)^{\frac{1}{2}}\|f\|_2.
\]
Thus Fubini's theorem applies locally,
and the formula for $Vf$ holds for
a.e. $x>0$.

Finally, for $f\in L^2(0,\infty)$, the
Cauchy--Schwarz inequality shows that
\[
\int_a^\infty e^{-u}\frac{|f(u)|}{u}\,du
\le \|f\|_2
\left(\int_a^\infty\frac{e^{-2u}}{u^2}\,du
\right)^{1/2}<\infty
\]
for every $a>0$. Hence $h=Vf$ is locally
absolutely continuous on $(0,\infty)$.
Differentiating \eqref{e20}  gives
\[
h'(x)=e^x\int_x^\infty e^{-u}
\frac{f(u)}{u}\,du-\frac{f(x)}{x}
=h(x)-\frac{f(x)}{x}
\]
for a.e. $x>0$, proving \eqref{e21}.
 \end{proof}

The semigroup representation of $\mathcal C^*$ yields its
pointwise integral formula on $H^2(\mathbb D)$ (see
\cite[(1.2)]{GallardoPartingtonRoss})
\begin{equation*}
(\mathcal C^*f)(z)=\frac{1}{1-z}\int_z^1 f(\zeta)\,d\zeta,
\qquad z\in\mathbb D.
\end{equation*}
 
Fix $\delta>0$ and let $B$ be a finite Blaschke product.
Let $\alpha_1,\ldots,\alpha_d\in\mathbb C_+$ be
distinct, with respective multiplicities
$n_1,\ldots,n_d\geq1$ and $d\in \mathbb{N}_+$.
Up to unimodular constants, set
\begin{equation}\label{e22}
b_\alpha(s):=\frac{s-\alpha}{s+\overline\alpha},
\qquad
B^+(s):=\prod_{j=1}^d b_{\alpha_j}(s)^{n_j},
\qquad
B:=B^+\circ q.
\end{equation}
Thus $B$ has zeros $b(\alpha_j)$ of multiplicity
$n_j$, and its degree is $N:=\sum_{j=1}^d n_j$. The Laplace transform allows us to represent the
finite-dimensional extension $K_{B\Theta_\delta}$ explicitly. Motivated
by the exponential--polynomial structure of
$\mathcal L^{-1}K_{B^+}$, define
\begin{equation}\label{e23}
E_{\delta,B}:=\bigvee\left\{
(x-\delta)^k e^{-\overline\alpha_j(x-\delta)}
\chi_{[\delta,\infty)}(x):
1\leq j\leq d,\ 0\leq k<n_j
\right\}.
\end{equation}
This is an $N$-dimensional subspace of
$L^2(0,\infty)$ by Lemma \ref{r12}. The model-space identity
$K_{uv}=K_u\oplus uK_v$, together with the
Paley--Wiener theorem and partial-fraction
decomposition, now yields
\begin{equation}\label{e24}
\mathcal J K_{B\Theta_\delta}
=L^2(0,\delta)\oplus E_{\delta,B}, \qquad \delta>0.
\end{equation}
The sum is orthogonal by support, while the
exponential--polynomial blocks for distinct
$\alpha_j$ are linearly independent but need
not be orthogonal.  Since $L^2(0,\delta)$ is invariant under $V$ and
every $Q_t$, the answers to the main questions concerning
$K_{B\Theta_\delta}$ are determined by \emph{their
action on the finite-dimensional space $E_{\delta,B}$.}

\section{Nearly invariant extensions and minimal subspaces}\label{s4}
\subsection{Minimal nearly invariant subspaces for the isometric shift}
The decomposition \eqref{e24} raises the question
of whether $\mathcal J K_{B\Theta_\delta}$ can be
generated by finitely many tail vectors through
the near-invariance condition. For $\delta>0$,
$\alpha\in\mathbb C_+$ and $n\in \mathbb{N}_+$, define
\[
f_{\delta,\alpha,n}(x):=
(x-\delta)^{n-1}e^{-\overline\alpha(x-\delta)}
\chi_{[\delta,\infty)}(x).
\]
The assumption $\delta>0$ is essential when $n>1$.
Indeed, $f_{0,\alpha,n}$ does not belong to
$\operatorname{ran}P_t$ for any $t>0$, so its
one-dimensional span is already nearly invariant,
but is strictly smaller than
$\mathcal L^{-1}K_{b_\alpha^n}$.

Modifying the finite-time orbit method of
\cite{LP22,LiangPartingtonCyclic}, we now
characterize the smallest closed nearly invariant
subspaces generated by these vectors, providing
a basis for the subsequent Ces\`aro defect calculations.

\begin{thm}\label{r5}
Let $\delta>0$, $\alpha\in\mathbb C_+$ and $n\in \mathbb{N}_+$.
The smallest closed subspace containing $f_{\delta,\alpha,n}$
and nearly invariant for $\{P_t\}_{t\ge0}$ is
$\mathcal L^{-1}K_{e^{-\delta s}b_\alpha^n}$.
If $d\in \mathbb{N}_+$ and $B^+=\prod_{j=1}^d b_{\alpha_j}^{n_j}$, where the
$\alpha_j$ are distinct, the smallest such subspace containing
$f_{\delta,\alpha_j,n_j}$ for  $1\le j \le d$ is
$\mathcal J K_{B\Theta_\delta}$.
\end{thm}

\begin{proof}

For $0\le u\le\delta$, we have
$P_{\delta-u}f_{u,\alpha,n}=f_{\delta,\alpha,n}$.
Hence every closed nearly invariant subspace containing
$f_{\delta,\alpha,n}$ also contains
\[
H:=\bigvee\{f_{u,\alpha,n}:0\le u\le\delta\}.
\]
We determine $H$ by computing its orthogonal complement.

Let $h\in H^\perp$. Then
$\langle f_{u,\alpha,n},h\rangle=0$
for $0\le u\le\delta$.
Repeated differentiation of the integral representation gives
\[
0=\left(\frac{\mathrm d}{\mathrm du}-\overline\alpha\right)^n
\langle f_{u,\alpha,n},h\rangle
=(-1)^n(n-1)!\overline{h(u)},
\qquad\text{for a.e. }u\in(0,\delta),
\]
implying that $h=0$ a.e. on $(0,\delta)$.

For $0\le u\le\delta$, expanding
$(x-u)^{n-1}=((x-\delta)+(\delta-u))^{n-1}$
and using
$e^{\overline\alpha(\delta-u)}
e^{-\overline\alpha(x-u)}
=e^{-\overline\alpha(x-\delta)}$, we obtain
\[
0=e^{\overline\alpha(\delta-u)}
\langle f_{u,\alpha,n},h\rangle=\sum_{k=0}^{n-1}\binom{n-1}{k}
(\delta-u)^{n-1-k}
\int_\delta^\infty
(x-\delta)^k e^{-\overline\alpha(x-\delta)}
\overline{h(x)}\,\mathrm dx.
\]
Since this polynomial vanishes for all
$u\in[0,\delta]$, its coefficients vanish.
Conversely, these conditions, together with
$h=0$ a.e. on $(0,\delta)$, imply
$\langle f_{u,\alpha,n},h\rangle=0$
for all $u\in[0,\delta]$. Therefore
\[
H=L^2(0,\delta)\oplus
\bigvee\big\{
(x-\delta)^k e^{-\overline\alpha(x-\delta)}
\chi_{[\delta,\infty)}(x):0\le k<n
\big\}=\mathcal L^{-1}K_{e^{-\delta s}b_\alpha^n},
\]
where the last equality follows from
\eqref{e23}--\eqref{e24} with $B^+=b_\alpha^n$.

To verify near invariance, use \eqref{mus} and recall that for any inner function
$\omega$ on $\mathbb C_+$, multiplication by $e^{-ts}$ is
isometric on $H^2(\mathbb C_+)$. If
$e^{-ts}F\in K_\omega$, then for every
$G\in H^2(\mathbb C_+)$,
\[
\langle F,\omega G\rangle
=\langle e^{-ts}F,e^{-ts}\omega G\rangle=0,
\]
so $F\in K_\omega$. Hence
$\mathcal L^{-1}K_\omega$ is nearly invariant for
$\{P_t\}_{t\ge0}$. In particular, $H$ is nearly invariant. \vspace{1mm}

For the second assertion, let $d\in \mathbb{N}_+$.
By the first part, every closed nearly invariant
subspace containing $f_{\delta,\alpha_j,n_j}$
for $1\le j\le d$ must contain
\[
 \bigvee\{f_{u,\alpha_j,n_j}:
0\le u\le\delta,\ 1\le j\le d\}
 =L^2(0,\delta)\oplus E_{\delta,B}
=\mathcal J K_{B\Theta_\delta},
\]
where the first equality follows by applying
the preceding orthogonal-complement calculation
to each $\alpha_j$, and the second follows
from \eqref{e24}. Meanwhile, 
$\mathcal J K_{B\Theta_\delta}
=\mathcal L^{-1}K_{e^{-\delta s}B^+}$
is itself nearly invariant for
$\{P_t\}_{t\ge0}$, as established above.
Hence $\mathcal J K_{B\Theta_\delta}$
is the smallest closed nearly invariant
subspace containing the prescribed vectors.
\end{proof}

\begin{rem}
The smallest closed \emph{nearly invariant and invariant
subspaces} generated by $f_{\delta,\alpha,n}$ under
$\{P_t\}_{t\geq0}$ are not comparable. Indeed,
Theorem~\ref{r5} gives
\[
M_{\mathrm{near}}:=\mathcal L^{-1}K_{e^{-\delta s}b_\alpha^n}\\
=L^2(0,\delta)\oplus
\bigvee\left\{
(x-\delta)^k e^{-\overline\alpha(x-\delta)}
\chi_{[\delta,\infty)}(x):0\le k<n
\right\}.
\]
On the other hand, since
\((\mathcal Lf_{\delta,\alpha,n})(s)
=\frac{(n-1)!e^{-\delta s}}
{(s+\overline\alpha)^n}
\) and  $(s+\overline\alpha)^{-n}$ is outer in
$H^2(\mathbb C_+)$, it yields that
\[
M_{\mathrm{inv}}
:=\bigvee\{P_t f_{\delta,\alpha,n}:t\ge0\}
=L^2(\delta,\infty).
\] Consequently, $M_{\mathrm{near}}$ contains
$L^2(0,\delta)$, which is orthogonal to
$M_{\mathrm{inv}}$, whereas
$M_{\mathrm{near}}\cap M_{\mathrm{inv}}$
is $n$-dimensional. Since $M_{\mathrm{inv}}$
is infinite dimensional, neither space contains
the other.
\end{rem}

In Theorem~\ref{r5}, the model space
$\mathcal J K_{B\Theta_\delta}$ is obtained as the
smallest closed nearly invariant subspace containing
one exponential--polynomial function for each
distinct zero of $B$. We now show that the same
space can be obtained from a single function.

\begin{cor}\label{r6}
Let $d\in \mathbb{N}_+$ and let $\alpha_1,\ldots,\alpha_d$ be distinct points
in $\mathbb C_+$, with positive integers $n_1,\ldots,n_d$. Set
\[
Q(s):=\prod_{j=1}^d(s+\overline\alpha_j)^{n_j},
\qquad g:=\mathcal L^{-1}(Q^{-1}).
\]
For every $\delta>0$, the smallest closed subspace
containing $P_\delta g$ and nearly invariant for
$\{P_t\}_{t\ge0}$ is $\mathcal J K_{B\Theta_\delta}$.
\end{cor}

\begin{proof}
Put $N:=\sum_{j=1}^d n_j$. The convolution formula for the inverse Laplace
transform gives
\[
g(x)=\frac{x^{N-1}}{(N-1)!}+O(x^N),
\qquad x\downarrow0.
\]
In particular, $g^{(k)}(0)=0$ for $0\le k<N-1$.
As $g$ is an exponential polynomial, all its
derivatives belong to $L^2(0,\infty)$. Integration
by parts therefore yields
\[
\mathcal L(g^{(k)})(s)=s^k \mathcal Lg(s)=\frac{s^k}{Q(s)},
\qquad 0\le k<N.
\]

The rational functions $s^k/Q(s)$,
$0\le k<N$, form a basis of the
$N$-dimensional space of proper rational
functions with poles among
$-\overline{\alpha_j}$ of orders at most $n_j$.
By partial-fraction decomposition, the same
space is spanned by
\[
\left\{
\frac{1}{(s+\overline{\alpha_j})^\ell}:
1\le j\le d,\ 1\le\ell\le n_j
\right\}.
\]
Since
\begin{align}\label{Lxs}
\mathcal L^{-1}
\left(\frac{1}{(s+\overline{\alpha_j})^\ell}\right)(x)
=
\frac{x^{\ell-1}}{(\ell-1)!}
e^{-\overline{\alpha_j}x}, 
\end{align}
applying the inverse Laplace transform
and then $P_u$ gives
\[
\bigvee_{k=0}^{N-1}P_ug^{(k)}
=
\bigvee\left\{
f_{u,\alpha_j,\ell}:
1\le j\le d,\ 1\le\ell\le n_j
\right\},
\qquad u\ge0.
\]
 
Now let
\(H:=\bigvee\{P_u g:0\le u\le\delta\}.
\) The vanishing conditions on the derivatives of
$g$ at $0$ ensure that the map
$u\mapsto P_u g$ is $(N-1)$ times continuously
differentiable in $L^2(0,\infty)$, with
\[
\frac{\mathrm d^k}{\mathrm du^k}P_u g
=(-1)^kP_u g^{(k)},
\qquad 0\le k<N.
\] 
Consequently, $P_u g^{(k)}\in H$ for
$0\le u\le\delta$ and $0\le k<N$,
where one-sided derivatives are used at the
endpoints. Thus
\[
f_{u,\alpha_j,n_j}\in H,
\qquad 0\le u\le\delta,\quad 1\le j\le d.
\]
The proof of Theorem~\ref{r5} now gives
\[
\mathcal J K_{B\Theta_\delta}
=\bigvee\left\{
f_{u,\alpha_j,n_j}:
0\le u\le\delta,\ 1\le j\le d
\right\}
\subseteq H.
\]

Conversely, $P_\delta g\in E_{\delta,B}
\subseteq\mathcal J K_{B\Theta_\delta}$.
Since $\mathcal J K_{B\Theta_\delta}$ is nearly
invariant for $\{P_t\}_{t\ge0}$ and
$P_{\delta-u}P_u g=P_\delta g$, we have
$P_u g\in\mathcal J K_{B\Theta_\delta}$
for every $0\le u\le\delta$. Hence
\( H=\mathcal J K_{B\Theta_\delta}.
\)

Finally, every closed nearly invariant subspace
containing $P_\delta g$ must contain $P_u g$
for $0\le u\le\delta$, and therefore contains
$H$. Since $H$ is itself nearly invariant,
the asserted minimality follows.
\end{proof}

The preceding isometric-shift statements can be viewed as special cases of the rational
orbit theorem in \cite[Theorem~4.10]{LiangPartingtonCyclic}. Indeed,
for $Q$ as in Corollary~\ref{r6}, put $G(s)=(s+1)^N/Q(s)$.
Both $G$ and $G^{-1}$ are bounded analytic on $\mathbb C_+$, and on
the imaginary axis $G/\overline G=B^+/b^N$. Consequently,
$G K_{e^{-\delta s}b^N}=K_{e^{-\delta s}B^+}$: this follows from
$T_{\overline{e^{-\delta s}B^+}}M_G
=T_{\overline G}T_{\overline{e^{-\delta s}b^N}}$
and invertibility of $T_{\overline G}$. Here $M_G$ denotes multiplication
by $G$, and $T_\psi f=P_+(\psi f)$ is the Toeplitz operator on
$H^2(\mathbb C_+)$ with bounded symbol $\psi$, where $P_+$ is the Hardy
projection.  

\subsection{Nearly invariant extensions for weighted translations}\label{e26}
%%%%%%%%%%%%%%%%%%%%%%%%%%%%%%%%%%%%%%%%%%%%%%%%%%%%%%%%%%%%%%%%%%%%%%%%%%%%

The construction for $\{P_t\}_{t\ge0}$ suggests a parallel
question for the weighted translations $\{S_t\}_{t\ge0}$.
Since the inverse translations are unbounded, the near-invariance
criterion must take account of their domains. We begin with the
support space $\mathcal N_a$ in \eqref{e7}.  The starting point is
the orthogonal decomposition
\[
L_w^2
=
\mathcal N_a
\oplus
L^2([a,\infty),w(y)\,dy),
\]
together with the  backward translation.
For $t>0$, define $ (\mathcal B_tf)(y):=f(y+t)
$  on its natural domain
\[
D(\mathcal B_t)
:=\ran S_t=
\left\{f\in L_w^2:\int_{\mathbb R}|f(x)|^2w(x-t)\,dx<\infty\right\},
\] and $\mathcal B_0=I$ on $L_w^2$.
Unlike $S_t$, $\mathcal B_t$ is generally unbounded on
$L_w^2$. Indeed,
\[
\frac{w(x-t)}{w(x)}
=
\exp\{2e^x(1-e^{-t})\}\longrightarrow\infty,
\qquad x\to+\infty.
\]
The operator $\mathcal B_t$ is closed and densely defined. Indeed,
$C_c^\infty(\mathbb R)\subseteq\operatorname{ran}S_t$, and if $f_n\to f$
and $\mathcal B_tf_n\to h$ in $L_w^2$ as $n\to\infty$, then
$f_n=S_t\mathcal B_tf_n\to S_th$ as $n\to\infty$. Hence $f=S_th\in D(\mathcal B_t)$
and $\mathcal B_tf=h$. Its restriction to $\mathcal N_a$ is well behaved:
$\mathcal B_t\mathcal N_a\subseteq\mathcal N_a$, and
\begin{equation}\label{e27}
\|\mathcal B_th\|_{L_w^2}
\leq
e^{\delta(1-e^{-t})}\|h\|_{L_w^2},
\qquad h\in\mathcal N_a.
\end{equation}

Near invariance of \(\mathcal N_a\oplus E\) requires control only of the part of the backward translation that remains on the complementary half-line \([a,\infty)\).  If $f\in L^2([a,\infty),w(y)\,dy)\cap D(\mathcal B_t)$,  then
\begin{equation} \label{e28}
\mathcal B_tf
=\chi_{[a-t,a)}f(\cdot+t)+\mathcal B_{a,t}f.
\end{equation}
The first term belongs automatically to $\mathcal N_a$; hence the
second term is precisely the part that must remain in the tail $E$. This motivates the truncated backward translation
\begin{equation}\label{e29}
(\mathcal B_{a,t}f)(y):=\chi_{[a,\infty)}(y)f(y+t),
\qquad y\in\mathbb R.
\end{equation}
Its domain is
$D(\mathcal B_{a,t}):=L^2([a,\infty),w(y)\,dy)
\cap D(\mathcal B_t)$, where functions on the half-line are
understood to be extended by zero to $\mathbb R$. Besides,
$\mathcal B_{a,0}=I$ on $L^2([a,\infty),w(y)\,dy)$.  

Using $\mathcal B_t$ and the truncated operator in \eqref{e29}, the near invariance of $\mathcal N_a\oplus E$  is determined exactly by the tail on which backward translation is admissible.

\begin{thm}\label{r7}
Let $a\in \mathbb{R}$ and $E$ be a closed subspace of
$L^2([a,\infty),w(y)\,dy)$. Then
\begin{equation} \label{e30}
\mathcal N_{a,E}:=\mathcal N_a\oplus E\subset L_w^2
\end{equation}
is nearly invariant for $\{S_t\}_{t\geq0}$
(equivalently, nearly invariant for $\{T_t\}_{t\geq0}$)
if and only if
\begin{equation}\label{e31}
\mathcal B_{a,t}\bigl(E\cap D(\mathcal B_t)\bigr)
\subseteq E,
\qquad t>0.
\end{equation}
 \end{thm}

\begin{proof}  Suppose first that \eqref{e31} holds and
$S_tg=h+e\in\mathcal N_a\oplus E$, where
$h\in\mathcal N_a$ and $e\in E$. Since
\[
\mathcal B_tS_tg=g\in L_w^2,
\]
we have $S_tg\in D(\mathcal B_t)$. Moreover,
$h\in\mathcal N_a\subset D(\mathcal B_t)$, and hence
\[
e=S_tg-h\in E\cap D(\mathcal B_t).
\]
Therefore,
\[
g=\mathcal B_tS_tg
=\mathcal B_th+\mathcal B_te.
\]
By \eqref{e27},
$\mathcal B_th\in\mathcal N_a$. Moreover,
\eqref{e28} gives
\begin{align}
\mathcal B_te
=
\chi_{[a-t,a)}e(\cdot+t)+\mathcal B_{a,t}e,\label{e32}
\end{align}
whose first term belongs to $\mathcal N_a$, while the second belongs
to $E$ by \eqref{e31}. Thus
$g\in\mathcal N_a\oplus E$, proving the near invariance.

Conversely, assume $\mathcal N_{a,E}$ is nearly invariant for $\{S_t\}_{t\geq0}$. Let
$e\in E\cap D(\mathcal B_t)$, then
$\mathcal B_te\in L_w^2$ and
\[
S_t\mathcal B_te=e\in E\subseteq\mathcal N_{a,E}.
\]
Near invariance therefore implies $  \mathcal B_te\in\mathcal N_{a,E}$. Together with  \eqref{e32}, this implies $\mathcal B_{a,t}e\in E$, since
$\mathcal B_{a,t}e$ is supported on $[a,\infty)$.  
  So \eqref{e31} holds. 
\end{proof}
We include an example to show it may happen that
$ E \cap D(\mathcal B_t)=\{0\},$ implying \eqref{e31} is trivial.
\begin{exm}\label{r8}
For $a\in \mathbb{R}$, take the function
\[
g_a(y) =
\frac{e^{e^y}}{1+(y-a)^2}
\chi_{[a,\infty)}(y)
\]
and define $E_a:=\mathbb C g_a$. Since
 \[
\|g_a\|_{L_w^2}^2
=
e^2\int_a^\infty
\frac{dy}{(1+(y-a)^2)^2}
=
e^2\int_0^\infty\frac{dx}{(1+x^2)^2}
=
\frac{\pi e^2}{4}<\infty,
\]
hence $g_a\in L_w^2$, and therefore 
$ E_a\subset L^2([a,\infty),w(y)\,dy).$ 

On the other hand, for $t>0$, the translated function $g_a(\cdot+t)$ satisfies
\[
\int_{\mathbb R}|g_a(y+t)|^2w(y)\,dy
\ge e^2\int_a^\infty
\frac{\exp\{2(e^t-1)e^y\}}
     {(1+(y+t-a)^2)^2}\,dy
=\infty.
\]
Thus $g_a(\cdot+t)\notin L_w^2$, and hence
$g_a\notin D(\mathcal B_t)$. Consequently, $E_a\cap D(\mathcal B_t)=\{0\}$ for every $t>0$, and  Theorem \ref{r7} shows that $\mathcal N_a\oplus E_a$ is nearly invariant for $\{S_t\}_{t\geq0}$.
\end{exm}
The preceding example shows that near invariance alone may impose
no effective restriction on  the tail space $E$ when
$E\cap D(\mathcal B_t)=\{0\}$. To rule out such vacuous extensions, we say that $\mathcal N_a\oplus E$ has a
\emph{backward-admissible tail} if
\begin{align}
E\subseteq D_\infty(\mathcal B)
:= \bigcap_{t>0}D(\mathcal B_t).\label{e33}
\end{align}   Theorem~\ref{r7} immediately yields the following corollary.
\begin{cor}\label{r9} Let $a\in \mathbb{R}$ and $E$ be a closed subspace of
$L^2([a,\infty),w(y)\,dy)$ with a backward-admissible tail, i.e., \eqref{e33} holds. Then
$\mathcal N_a\oplus E$ 
is nearly invariant for $\{S_t\}_{t\geq0}$
(equivalently, nearly invariant for $\{T_t\}_{t\geq0}$)
if and only if
\begin{equation*}
\mathcal B_{a,t}E\subseteq E,
\qquad t>0.
\end{equation*}
\end{cor}
\medskip

Corollary~\ref{r9} suggests constructing nearly invariant
extensions of $\mathcal N_a$ by adjoining finite-dimensional
tail spaces that are invariant under truncated backward
translation. This parallels the finite Blaschke extensions
of model spaces. Indeed, by \eqref{e17},
$\mathcal N_a$ corresponds to the singular model space
$K_{e^{-\delta s}}$, where $\delta=e^a$, while
\[
K_{e^{-\delta s}b^n}
=K_{e^{-\delta s}}\oplus e^{-\delta s}K_{b^n}
\]
is a finite-dimensional extension of $K_{e^{-\delta s}}$.
Motivated by this structure, we consider analogous
extensions of $\mathcal N_a$ generated by
\emph{exponential--polynomial tails}, whose backward
translates are finite linear combinations of the
same functions modulo $\mathcal N_a$.

For $a\in\mathbb R$, $\lambda\in\mathbb C$ and
$k\in\mathbb N_0$, define
\begin{equation}\label{e34}
r_{a,\lambda,k}(y):=
\frac{(y-a)^k}{k!}e^{\lambda(y-a)}
\chi_{[a,\infty)}(y).
\end{equation}
The following lemma ensures that these functions
belong to $L_w^2$.

\begin{lem}\label{r10}
Let $a\in\mathbb R$, $\lambda\in\mathbb C$ and $k\in\mathbb N_0$. Then
$r_{a,\lambda, k}\in L_w^2$.
\end{lem}

\begin{proof}
With $\delta=e^a$ and $\alpha=\Re\lambda$, the substitution $x=e^y$
gives
\[
\|r_{a,\lambda,k}\|_{L_w^2}^2
=\frac{e^{2-2\alpha a}}{(k!)^2}
\int_\delta^\infty (\log x-a)^{2k}x^{2\alpha-1}e^{-2x}\,dx<\infty.
\]
The integral converges at its positive lower endpoint and at infinity
because the exponential factor dominates every power of $x$ and
$\log x$.
\end{proof}
 
We next show that  every $r_{a,\lambda,k}$ belongs to the domain of each backward translation $\mathcal{B}_{t}.$
\begin{lem}\label{r11}
    Let $a\in\mathbb R$, $\lambda\in\mathbb C$,
$k\in\mathbb N_0$, and $t>0$,  then
$r_{a,\lambda,k}\in D(\mathcal B_t)$ and
\begin{equation}\label{e35}
\mathcal B_{a,t}r_{a,\lambda,k}
=
e^{\lambda t}
\sum_{j=0}^k
\frac{t^{k-j}}{(k-j)!}\,r_{a,\lambda,j}.
\end{equation}
\end{lem}

\begin{proof}
For $y\ge a$, the binomial formula gives the right-hand side of
\eqref{e35}. It belongs to $L_w^2$ by Lemma~\ref{r10}. On the remaining
interval $[a-t,a)$ the backward translate is bounded and the weight is
bounded, so it is square integrable there as well. It vanishes on
$(-\infty,a-t)$. Thus $r_{a,\lambda,k}\in D(\mathcal B_t)$, and
restriction to $[a,\infty)$ proves \eqref{e35}.
\end{proof}

To collect finite-dimensional Jordan blocks with distinct exponential
rates, let \(
\Lambda:=\{\lambda_1,\ldots,\lambda_r\}\subset\mathbb C
\) be a finite set of distinct points and let
$\bm m:=(m_1,\ldots,m_r)\in\mathbb N_+^r$. Using \eqref{e34}, define the associated tail subspace
\begin{equation*}
E_{a,\Lambda,\bm m}
:=\bigvee \{r_{a,\lambda_\nu,k}: \;1\leq\nu\leq r,\; 0\leq k<m_\nu \}.
\end{equation*}

\begin{lem}\label{r12}
Let $a\in\mathbb R$,  $r\in\mathbb N_+$,
$\Lambda=\{\lambda_1,\ldots,\lambda_r\}\subset\mathbb C$
consist of distinct points, and  
$\bm m=(m_1,\ldots,m_r)\in\mathbb N_+^r$. Then
\begin{equation}\label{e38}
\dim E_{a,\Lambda,\bm m}
=
\sum_{\nu=1}^r m_\nu.
\end{equation}
\end{lem}

\begin{proof}
It suffices to prove that the functions
 $r_{a,\lambda_\nu,k},\;
1\leq\nu\leq r,\; 0\leq k<m_\nu,$ 
are linearly independent. Suppose that
\begin{align}\sum_{\nu=1}^r\sum_{k=0}^{m_\nu-1}c_{\nu,k}r_{a,\lambda_\nu,k}= \sum_{\nu=1}^r\sum_{k=0}^{m_\nu-1}c_{\nu,k} \frac{(y-a)^k}{k!}
e^{\lambda_\nu (y-a)}\chi_{[a,\infty)}(y)=0.\label{e39}\end{align}
For $x=y-a\geq0$, set
\begin{align*}
p_\nu(x):=
\sum_{k=0}^{m_\nu-1}
c_{\nu,k}\frac{x^k}{k!},
\qquad 1\leq\nu\leq r.
\end{align*}
The equality in \eqref{e39} is initially an equality almost
everywhere. Since the finite exponential--polynomial sum is
continuous on $[0,\infty)$, it vanishes everywhere there. Hence
\begin{align}
\sum_{\nu=1}^r
p_\nu(x)e^{\lambda_\nu x}=0,
\qquad x\geq0.\label{e41}
\end{align}
Fix $\mu\in\{1,\ldots,r\}$ and apply the differential operator
$\prod_{\nu\ne\mu}(\frac{\mathrm d}{\mathrm d x}-\lambda_\nu)^{m_\nu}$
to \eqref{e41}. Every term with index $\nu\ne\mu$ vanishes,
because the polynomial $p_\nu$ has degree less than $m_\nu$.
For the remaining term, we obtain
\begin{equation}\label{e42}
\prod_{\nu\ne\mu}
\left(\frac{\mathrm d}{\mathrm d x}+\lambda_\mu-\lambda_\nu\right)^{m_\nu}p_\mu=0.
\end{equation}
For $\nu\ne\mu$, the constant $\lambda_\mu-\lambda_\nu$ is
nonzero. Differentiation is nilpotent on the space of polynomials
of degree less than $m_\mu$, so each factor in \eqref{e42}
is invertible on that space. Hence $p_\mu=0$. Since $\mu$ was
arbitrary, all $p_\nu$ vanish, and therefore
\[
c_{\nu,k}=0,
\qquad
1\leq\nu\leq r,\quad 0\leq k\leq m_\nu-1.
\]
Thus the dimension formula \eqref{e38}  follows.
\end{proof}

For distinct parameters
$\Lambda=\{\lambda_1,\ldots,\lambda_r\}$ and multiplicities
$\bm m=(m_1,\ldots,m_r)$, set
\begin{equation}\label{e43}
\mathcal N_{a,\Lambda,\bm m} :=
\mathcal N_a\oplus E_{a,\Lambda,\bm m},
\end{equation} where the sum is orthogonal because $\mathcal N_a$ and
$E_{a,\Lambda,\bm m}$ are supported, up to a null set, on
$(-\infty,a)$ and $[a,\infty)$, respectively.

\begin{thm}\label{r13} Let $a\in\mathbb R$,  $r\in\mathbb N_+$, $\Lambda=\{\lambda_1,\ldots,\lambda_r\}\subset\mathbb C$
consist of distinct points, and  
$\bm m=(m_1,\ldots,m_r)\in\mathbb N_+^r$. Then $\mathcal N_{a,\Lambda,\bm m}$ in \eqref{e43}  is nearly invariant for $\{S_t\}_{t\geq0}$. Moreover,
\[
\dim(\mathcal N_{a,\Lambda,\bm m}\ominus \mathcal N_a)
=\sum_{\nu=1}^r m_\nu.
\]
\end{thm}

\begin{proof}
Lemma~\ref{r11} gives  \[\mathcal B_{a,t}r_{a,\lambda_\nu,k}
=e^{\lambda_\nu t} \sum_{j=0}^k \frac{t^{k-j}}{(k-j)!}\,
r_{a,\lambda_\nu,j}, \qquad 1\leq\nu\leq r,\; 0\leq k<m_\nu.\]
Since $j\leq k<m_\nu$, every term on the right-hand side belongs
to $E_{a,\Lambda,\bm m}$. Consequently,
\begin{equation*}
\mathcal B_{a,t}E_{a,\Lambda,\bm m}
\subseteq E_{a,\Lambda,\bm m},
\qquad t>0.
\end{equation*}
Together with Theorem~\ref{r7}, this implies
that $\mathcal N_{a,\Lambda,\bm m}$ is nearly invariant for $\{S_t\}_{t\geq0}$.
Since the sum in \eqref{e43} is orthogonal,
Lemma~\ref{r12} yields
\[
\dim(\mathcal N_{a,\Lambda,\bm m}\ominus \mathcal N_a)
=
\dim E_{a,\Lambda,\bm m}
=
\sum_{\nu=1}^r m_\nu. \qedhere
\] \end{proof}

The construction \eqref{e43} also has an explicit
$H^2(\mathbb C_+)$ interpretation. By the unitary map $\mathcal F$ in \eqref{e14} and \eqref{e17}, it follows that 
\[
U\mathcal F^{-1}\mathcal N_a=K_{e^{-\delta s}}.
\]
 For $r_{a,\lambda,k}$ in \eqref{e34}, set
\begin{equation}\label{e44}
\psi_{\delta,\lambda,k}(s)
:= \frac{1}{k!}\int_\delta^\infty
\left(\log\frac{x}{\delta}\right)^k
x^{\lambda-\frac12}e^{-(s+1)x}\,dx.
\end{equation}
A direct calculation gives
$U\mathcal F^{-1} r_{a,\lambda,k}
=e\delta^{-\lambda}\psi_{\delta,\lambda,k}$. Hence set
\begin{equation*}
\mathcal E_{\delta,\Lambda,\bm m} :=
\bigvee\left\{
\psi_{\delta,\lambda_\nu,k}:
1\leq\nu\leq r,\ 0\leq k<m_\nu
\right\}\subset H^2(\mathbb{C}_+).
\end{equation*}
Then
\begin{equation}\label{e45}
U\mathcal F^{-1}\mathcal N_{a,\Lambda,\bm m}
= K_{e^{-\delta s}}\oplus
\mathcal E_{\delta,\Lambda,\bm m}.
\end{equation}
Thus $\mathcal N_{a,\Lambda,\bm m}$ corresponds on
$H^2(\mathbb C_+)$ to an explicit finite-dimensional extension of the
singular model space $K_{e^{-\delta s}}$. \emph{In general this extension
need not itself be a model space.} The precise model-space case will be
characterized below.

Although the exponential--polynomial extensions
$\mathcal N_{a,\Lambda,\bm m}$ in \eqref{e43} are nearly invariant for
$\{S_t\}_{t\ge0}$, they are invariant under
neither $\{S_t\}_{t\ge0}$ nor   $\{S_t^*\}_{t\ge0}$, as the following
two propositions show.

\begin{prop}\label{r14}
Let $a\in\mathbb R$,  $r\in\mathbb N_+$, 
$\Lambda=\{\lambda_1,\ldots,\lambda_r\}\subset\mathbb C$
consist of distinct points, and  
$\bm m=(m_1,\ldots,m_r)\in\mathbb N_+^r$. Then, for every $t>0$,
$S_t^*\mathcal N_{a,\Lambda,\bm m}
\nsubseteq
\mathcal N_{a,\Lambda,\bm m}.$
\end{prop}

\begin{proof}
Fix $t>0$ and choose any
$\lambda_\mu\in\Lambda$. Since
$r_{a,\lambda_\mu,0}\in E_{a,\Lambda,\bm m}$, it follows from
\eqref{e16} that, for $y\geq a$,
\begin{align*}
(S_t^*r_{a,\lambda_\mu,0})(y)
 =e^{-2e^y(e^t-1)}
r_{a,\lambda_\mu,0}(y+t)  =
e^{\lambda_\mu t}
e^{-2e^y(e^t-1)}
e^{\lambda_\mu(y-a)}.
\end{align*}
Suppose, on the contrary, that
$S_t^*r_{a,\lambda_\mu,0}
\in\mathcal N_{a,\Lambda,\bm m}$.
Restricting to  $[a,\infty)$ then gives  
\(S_t^*r_{a,\lambda_\mu,0} \in E_{a,\Lambda,\bm m},\) i.e., 
\begin{align}
e^{\lambda_\mu (t-a)}
e^{\lambda_\mu y} e^{-2e^y(e^t-1)} 
=\sum_{\nu=1}^r
p_\nu(y-a)e^{\lambda_\nu(y-a)},\qquad
\deg p_\nu<m_\nu \label{e46}
\end{align}
for some polynomials $p_\nu$. Since
\[
(D-\lambda_\nu)^{m_\nu}
\bigl(
p_\nu(y-a)e^{\lambda_\nu(y-a)}
\bigr)
= p_\nu^{(m_\nu)}(y-a)
e^{\lambda_\nu(y-a)}
=0\] for $D=\frac{\mathrm d}{\mathrm d y}$. Applying
$\prod_{\nu=1}^r(D-\lambda_\nu)^{m_\nu}$ to both sides of
\eqref{e46}, the right-hand side vanishes identically. On the
other hand, set $c=2(e^t-1)>0$. For every polynomial $P$,
\[
(D-\lambda)
\left(
e^{\lambda_\mu y}e^{-ce^y}P(e^y)
\right)
=
e^{\lambda_\mu y}e^{-ce^y}Q(e^y),
\]
where
$Q(x)=x\frac{\mathrm d P}{\mathrm d x}(x)+(\lambda_\mu-\lambda-cx)P(x).$ If $P$ has degree $d$ and leading coefficient $a_d\neq0$, then
$Q$ has degree $d+1$ and leading coefficient $-ca_d\neq0$.
Starting with $P\equiv1$ and iterating this observation yields
\[
\prod_{\nu=1}^r(D-\lambda_\nu)^{m_\nu}
\left(
e^{\lambda_\mu t}
e^{-2e^y(e^t-1)}
e^{\lambda_\mu(y-a)}
\right)
\not\equiv0,
\]
which contradicts \eqref{e46}. Thus
$ S_t^*r_{a,\lambda_\mu,0}
\notin \mathcal N_{a,\Lambda,\bm m},
$ implying the desired result. 
\end{proof}

\begin{prop}\label{r15}
 Let $a\in\mathbb R$,  $r\in \mathbb{N}_+$,  
$\Lambda=\{\lambda_1,\ldots,\lambda_r\}\subset\mathbb C$
consist of distinct points, and  
$\bm m=(m_1,\ldots,m_r)\in\mathbb N_+^r$. Then, for every $t>0$,
\(
S_t\mathcal N_{a,\Lambda,\bm m}
\nsubseteq
\mathcal N_{a,\Lambda,\bm m}.
\)
\end{prop}

\begin{proof}
Fix $t>0$ and choose a nonzero
$g\in L^2((a-t,a),w(y)\,dy)$. Then $g\in\mathcal N_a$, whereas
$S_tg$ is supported in $(a,a+t)$.  If $S_tg\in\mathcal N_{a,\Lambda,\bm m}$, then its
$\mathcal N_a$-component vanishes, so
$S_tg\in E_{a,\Lambda,\bm m}$. On $[a,\infty)$ every element of
$E_{a,\Lambda,\bm m}$ is an exponential polynomial
\[
F(y):=
\sum_{\nu=1}^r
p_\nu(y-a)e^{\lambda_\nu(y-a)},
\qquad
\deg p_\nu<m_\nu.
\]
Since $S_tg$ vanishes a.e. on $(a+t,\infty)$,
continuity implies that $F$ vanishes there identically. Analyticity
then gives $F\equiv0$, contradicting $S_tg\neq0$.
\end{proof}

We next identify the Jordan structure of the truncated backward
translations on $E_{a,\Lambda,\bm m}$.  For $a\in \mathbb{R}, \lambda\in\mathbb C$ and $m\in\mathbb N_+$, using \eqref{e34}, set
\begin{equation} 
E_{a,\lambda,m}:=
\bigvee \{r_{a,\lambda,k}: \; 0\leq k<m  \},\label{ealm}
\end{equation}
and consider the ordered basis
\(\mathscr B_{\lambda,m}:=
\bigl(
r_{a,\lambda,0},
r_{a,\lambda,1},
\ldots,
r_{a,\lambda,m-1}
\bigr).
\)
By \eqref{e35},  $\{\mathcal B_{a,t}|_{E_{a,\lambda,m}}\}_{t\geq0}$
forms a finite-dimensional semigroup.

Let $\mathcal G_{a,\lambda,m}$ denote its infinitesimal generator.
Since $E_{a,\lambda,m}$ is finite dimensional, it is given by
\[
\mathcal G_{a,\lambda,m}f
=
\left.
\frac{\mathrm d}{\mathrm d t}
\mathcal B_{a,t}f
\right|_{t=0},
\qquad
f\in E_{a,\lambda,m}.
\]
Differentiating \eqref{e35}  at $t=0$ gives
\begin{equation*}
\mathcal G_{a,\lambda,m}r_{a,\lambda,0}
=
\lambda r_{a,\lambda,0},
\qquad
\mathcal G_{a,\lambda,m}r_{a,\lambda,k}
=
\lambda r_{a,\lambda,k}
+r_{a,\lambda,k-1},
\quad 1\leq k<m.
\end{equation*}
Therefore, relative to the ordered basis
$\mathscr B_{\lambda,m}$,
\begin{equation*}
[\mathcal G_{a,\lambda,m}]_{\mathscr B_{\lambda,m}}
=
J_{\lambda,m},
\end{equation*}
where
\begin{equation*}
J_{\lambda,m}
=
\begin{pmatrix}
\lambda&1&0&\cdots&0\\
0&\lambda&1&\ddots&\vdots\\
\vdots&\ddots&\ddots&\ddots&0\\
0&\cdots&0&\lambda&1\\
0&\cdots&\cdots&0&\lambda
\end{pmatrix}.
\end{equation*}
Equivalently, since
\(
\mathcal B_{a,t}\big|_{E_{a,\lambda,m}}
=e^{t\mathcal G_{a,\lambda,m}},
\)
we obtain
\begin{equation*}
\left[
\mathcal B_{a,t}\big|_{E_{a,\lambda,m}}
\right]_{\mathscr B_{\lambda,m}}
=
e^{tJ_{\lambda,m}}.
\end{equation*}

For the general extension, Lemma~\ref{r12} gives the
algebraic direct sum
\begin{align}
E_{a,\Lambda,\bm m}
=
\dotplus_{\nu=1}^r
E_{a,\lambda_\nu,m_\nu}.\label{e47}
\end{align}
Let
\(
\mathscr B_{\Lambda,\bm m}
= \bigl(
\mathscr B_{\lambda_1,m_1},
\ldots,
\mathscr B_{\lambda_r,m_r}
\bigr)
\)
be the corresponding ordered basis. Since each
$E_{a,\lambda_\nu,m_\nu}$ is invariant under
$\mathcal B_{a,t}$, the generator
$\mathcal G_{a,\Lambda,\bm m}$ of
$\{\mathcal B_{a,t}|_{E_{a,\Lambda,\bm m}}\}_{t\geq0}$
satisfies
\begin{equation*} 
[\mathcal G_{a,\Lambda,\bm m}]_{\mathscr B_{\Lambda,\bm m}}
=
J_{\Lambda,\bm m}
:=
\bigoplus_{\nu=1}^r
J_{\lambda_\nu,m_\nu}.
\end{equation*}
Consequently,
\begin{equation*}
\left[
\mathcal B_{a,t}\big|_{E_{a,\Lambda,\bm m}}
\right]_{\mathscr B_{\Lambda,\bm m}}
=
e^{tJ_{\Lambda,\bm m}}
=
\bigoplus_{\nu=1}^r
e^{tJ_{\lambda_\nu,m_\nu}}.
\end{equation*}
In particular, $\sigma(\mathcal G_{a,\Lambda,\bm m})=\Lambda.$ 
Thus the prescribed Jordan data \((\lambda_\nu,m_\nu)\) are realized by the finite-dimensional backward-translation dynamics.

\subsection{Model-space rigidity for weighted-translation extensions}\label{s5}
Considering the important role of the model space, we now determine exactly \emph{which exponential--polynomial extensions
$\mathcal N_{a,\Lambda,\bm m}$ in \eqref{e43} become model spaces in Hardy space} on the right half-plane $H^2(\mathbb{C}_+)$ under the
unitary map $U\mathcal F^{-1}$. Within this family, the model-space condition imposes specific
restrictions on the exponential rates and their multiplicities.

\begin{thm}\label{r17} Let $a\in\mathbb R$, $r\in\mathbb{N}_+$,  
$\Lambda=\{\lambda_1,\ldots,\lambda_r\}\subset\mathbb C$
consist of distinct points, and  
$\bm m=(m_1,\ldots,m_r)\in\mathbb N_+^r$ such that
$N=\sum_{\nu=1}^r m_\nu$. Then $U\mathcal F^{-1}\mathcal N_{a,\Lambda,\bm m}$ is a model space in $H^2(\mathbb C_+)$ if and only if, after
a reordering of the parameters,
\begin{align}
r=N,\qquad
m_1=\cdots=m_N=1,\qquad
\lambda_j=j-\frac12,\quad 1\leq j\leq N.\label{e48}
\end{align}
In this case,
 $U\mathcal F^{-1}\mathcal N_{a,\Lambda,\bm1}=K_{e^{-\delta s}b(s)^N}.$
\end{thm}

\begin{proof}
\emph{Necessity.}
Assume that
$U\mathcal F^{-1}\mathcal N_{a,\Lambda,\bm m}=K_\Phi$
for some inner function $\Phi$ on $\mathbb C_+$. Since
$K_{e^{-\delta s}}\subseteq K_\Phi$,   the corresponding orthogonal complements satisfy
 $ \Phi H^2(\mathbb C_+)
\subseteq
e^{-\delta s}H^2(\mathbb C_+). $ Hence, up to a unimodular constant, $ \Phi=e^{-\delta s}B$ for some inner function $B$ on $\mathbb C_+$.  
Using the standard identity
$K_{e^{-\delta s}B}=K_{e^{-\delta s}}\oplus e^{-\delta s} K_B
$  and comparing with
\eqref{e45}, we obtain 
\begin{equation}\label{e49}
\mathcal E_{\delta,\Lambda,\bm m}
=
e^{-\delta s}K_B.
\end{equation}
By Lemma~\ref{r12},
$\dim\mathcal E_{\delta,\Lambda,\bm m}=N$, and hence
$\dim K_B=N$. Consequently, $B$ is a finite Blaschke product of
degree $N$.

We next use the Laplace-transform representation of the vectors
in $\mathcal E_{\delta,\Lambda,\bm m}$. For
$\lambda\in\mathbb C$ and $k\in\mathbb N_0$, the change of
variables $x=\delta+t$ in
\eqref{e44} gives
\begin{align*}
\psi_{\delta,\lambda,k}(s)
=
\frac{e^{-\delta(s+1)}}{k!}
\int_0^\infty
(\delta+t)^{\lambda-\frac12}
\left(
\log\left(1+\frac{t}{\delta}\right)
\right)^k
e^{-(s+1)t}\,dt \notag=
e^{-\delta}e^{-\delta s}
\mathcal L h_{\lambda,k}(s),
\end{align*}
where
\begin{align*}
h_{\lambda,k}(t):=
\frac{1}{k!}
(\delta+t)^{\lambda-\frac12}
\left(
\log\left(1+\frac{t}{\delta}\right)
\right)^k e^{-t}.
\end{align*}
It follows from \eqref{e49} that
\begin{equation}\label{e50}
\mathcal L^{-1} (K_B)
= \bigvee
\left\{
 h_{\lambda_\nu,k}:
1\leq\nu\leq r,\ 0\leq k<m_\nu
\right\}.
\end{equation}

Let $\alpha_1,\ldots,\alpha_q\in\mathbb C_+$ be the distinct
zeros of $B$, with multiplicities $d_1,\ldots,d_q$. Thus
$d_1+\cdots+d_q=N$, and
\[
K_B
=
\bigvee
\left\{
\frac{1}{(s+\overline{\alpha_\ell})^j}:
1\leq\ell\leq q,\ 1\leq j\leq d_\ell
\right\}.
\]
By \eqref{Lxs}, we have
\begin{equation}\label{e51}
\mathcal L^{-1}(K_B)=
\bigvee
\left\{
t^{j-1}e^{-\overline{\alpha_\ell}t}:\;
1\leq\ell\leq q,\ 1\leq j\leq d_\ell
\right\}.
\end{equation}
Every vector in the right-hand side of
\eqref{e51} extends to an entire function of $t$.
Taking $k=0$ in \eqref{e50} therefore shows, for every
$\nu$, that
$(\delta+t)^{\lambda_\nu-\frac12}e^{-t}$ 
has an entire continuation.  Since a noninteger power of $z+\delta$ is
multivalued near $-\delta$, single-valuedness
requires $\lambda_\nu-\frac12$ to be an integer.
Moreover, analyticity at $-\delta$ excludes
negative integer powers. Hence
$\lambda_\nu-\frac12\in\mathbb N_0$.
 Thus write
\[
\lambda_\nu=j_\nu+\frac12,
\qquad j_\nu\in\mathbb N_0.
\]
Since the $\lambda_\nu$ are distinct, so are the $j_\nu$.

We next exclude logarithmic multiplicities. If $m_\nu\geq2$, then
\eqref{e50} contains
\[
h_{\lambda_\nu,1}(t)
=
(\delta+t)^{j_\nu}
\log\!\left(1+\frac{t}{\delta}\right)e^{-t}.
\]
After one circuit around $t=-\delta$, the logarithm changes by
$2\pi i$, so this function cannot have a single-valued entire
continuation. This contradicts \eqref{e51}.
Therefore, $m_\nu=1$ for $ 1\leq\nu\leq r$ 
and hence $r=N$. Relabeling the distinct integers as
$j_1,\ldots,j_N$ gives
\begin{equation*}
\lambda_\nu=j_\nu+\frac12,
\qquad 1\leq\nu\leq N.
\end{equation*}
Hence \eqref{e50} yields
\begin{equation}\label{e52}
\mathcal L^{-1}(K_B)
=\bigvee
\left\{e^{-t}
(\delta+t)^{j_\nu}:1\leq\nu\leq N
\right\}.
\end{equation}

Comparing \eqref{e52} and \eqref{e51}, we see that
$e^{-\overline{\alpha_\ell}t}$ belongs to
$e^{-t}\mathcal P$ for a finite-dimensional polynomial space
$\mathcal P$. Thus $e^{(1-\overline{\alpha_\ell})t}$ is a
polynomial in $t$, which is possible only if
$\alpha_\ell=1$ for every $\ell$. Hence $B$ has the single
zero $1$ of multiplicity $N$.  Thus, up to a unimodular constant,
$B=b^N,$  where $b$ is given in \eqref{e8}. Consequently,
\[
K_B=K_{b(s)^N}
=
\bigvee
\left\{
\frac{1}{(1+s)^j}:1\leq j\leq N
\right\},
\]
and consequently
\[
\mathcal L^{-1}(K_B)
=
e^{-t}\mathcal P_{N-1},
\qquad
\mathcal P_{N-1}
=\bigvee\{1,t,\ldots,t^{N-1}\}.
\]
Together with \eqref{e52}, this gives
\[
\bigvee
\left\{
(\delta+t)^{j_\nu}:1\leq\nu\leq N
\right\}
=
\mathcal P_{N-1}.
\]
Since the integers $j_1,\ldots,j_N$ are distinct and the
polynomials $(\delta+t)^{j_\nu}$ have degrees $j_\nu$, the
equality above is possible only if
\[
\{j_1,\ldots,j_N\}
=
\{0,1,\ldots,N-1\}.
\]
After reordering, $j_\nu=\nu-1,$ and hence
$\lambda_\nu=\nu-\frac12, \; 1\leq\nu\leq N,$  
  proving the necessity.

\emph{Sufficiency.}
Conversely, suppose that \eqref{e48} holds, equivalently,
\[
\Lambda=
\left\{
j+\frac12:\;0\leq j\leq N-1
\right\}.
\]
For $0\leq j\leq N-1$, we have
\begin{equation}\label{e53}
\psi_{\delta,j+\frac12,0}(s)
= \int_\delta^\infty x^j e^{-(s+1)x}\,dx =
e^{-\delta(s+1)}
\sum_{k=0}^j
\frac{j!}{(j-k)!}
\frac{\delta^{j-k}}{(1+s)^{k+1}}.
\end{equation}
Since the coefficient of $(1+s)^{-(j+1)}$ in
\eqref{e53} is $j!\neq0$, these identities
give a triangular change of basis. Hence
\[
\mathcal E_{\delta,\Lambda,\bm1}
=
\bigvee
\left\{
\frac{e^{-\delta s}}{(1+s)^j}:1\leq j\leq N
\right\}
=
e^{-\delta s}K_{b(s)^N}.
\]
Using \eqref{e45}, we conclude that
\[
U\mathcal F^{-1}\mathcal N_{a,\Lambda,\bm1}
=
K_{e^{-\delta s}}
\oplus e^{-\delta s}K_{b(s)^N}=
K_{e^{-\delta s}b(s)^N},
\]
which indeed is a model space.
\end{proof}

  The rigidity of the model-space condition  is interpreted
 in terms of the Jordan dynamics as follows.

\begin{rem}\label{r18}
 In general, recall the  direct sum
\[
J_{\Lambda,\bm m}
:=
\bigoplus_{\nu=1}^r
J_{\lambda_\nu,m_\nu}
\]
  and a nontrivial Jordan block $J_{\lambda_\nu,m_\nu}$,
$m_\nu>1$, corresponds to the generalized eigenvectors
\[
(y-a)^k e^{\lambda_\nu(y-a)},
\qquad 1\leq k<m_\nu.
\]
After the logarithmic change of variables and the Laplace
transform, these generalized eigenvectors produce logarithmic
terms  of the form
\[
(\delta+t)^{\lambda_\nu-\frac12}
\left(
\log\left(1+\frac{t}{\delta}\right)
\right)^k e^{-t},\qquad 1\leq k<m_\nu.
\]
 Such functions cannot belong to the inverse-Laplace image of a
finite-dimensional model space, which only consists of
entire exponential polynomials. Thus the model-space condition
eliminates every nontrivial Jordan block.

Moreover, diagonalizability alone is not sufficient. The
eigenvalues are also rigidly determined: they must be exactly
$\frac12,\frac32,\ldots,N-\frac12.$  Hence, within the family $\mathcal N_{a,\Lambda,\bm m}$,
\emph{the model-space property is equivalent to the distinguished
Jordan configuration.} More precisely, after reordering the
parameters so that $\lambda_j=j-\frac12$, $1\leq j\leq N$,
one has
\[
J_{\Lambda,\bm 1}
=
\operatorname{diag}
\left(
\frac12,\frac32,\ldots,N-\frac12
\right).
\] This shows that the family
$K_{e^{-\delta s}b(s)^N}$
is precisely the model-space family within the class  of
exponential--polynomial extensions considered above.
\end{rem}
Transferring Theorem \ref{r17} back to the disk gives the following corollary. 
\begin{cor}\label{r19}
Let $a\in \mathbb{R}$ and $N\in\mathbb N_+$. Define
\begin{equation*}
\Lambda_N=
\left\{
\frac12,\frac32,\ldots,N-\frac12
\right\},
\qquad
\bm1=(1,\ldots,1).
\end{equation*}
Then  
 $\mathcal F^{-1}\mathcal N_{a,\Lambda_N,\bm1}
= K_{z^N\Theta_\delta}.
$ Conversely, among all spaces
$\mathcal N_{a,\Lambda,\bm m}$ with
$\sum_{\nu=1}^r m_\nu=N$, the space
$\mathcal F^{-1}\mathcal N_{a,\Lambda,\bm m}$ is a model space
in $H^2(\mathbb D)$ if and only if, up to reordering, $\Lambda=\Lambda_N$ and $\bm m=\bm 1$.
\end{cor}
We also note that the spaces $K_{z^N\Theta_\delta}$ are nearly invariant for $\{C_t\}_{t\geq0}$ and form the increasing
sequence 
\[
K_{\Theta_\delta}
\subsetneq K_{z\Theta_\delta}
\subsetneq K_{z^2\Theta_\delta}
\subsetneq\cdots,
\]
with
 $\dim\bigl(
K_{z^N\Theta_\delta}\ominus K_{\Theta_\delta}
\bigr)=N. $

\section{Almost invariance and Ces\`aro defects}
\subsection{Finite-dimensional Ces\`aro defects and single-generator extensions}
\label{s6}
 For every $f\in H^2(\mathbb D)$, equations \eqref{e11}--\eqref{e12} give
$(\mathcal Ff)(y)=e^{y/2}(\mathcal Jf)(e^y)/\sqrt{w(y)}$ for almost every $y\in\mathbb R$.
Setting $\mathfrak C:=\mathcal F\mathcal C\mathcal F^{-1}$, \eqref{e4}, \eqref{e15} and \eqref{e16} imply that
\begin{equation}\label{e56}
(\mathfrak Cg)(y)
=
e^{y/2+2e^y}
\int_y^\infty
e^{-s/2-2e^s}g(s)\,ds.
\end{equation}
In particular, $\mathfrak C\mathcal N_a\subseteq\mathcal N_a.$ 
Given a finite-dimensional subspace
$E\subset L^2([a,\infty),w(y)\,dy)$, recalling 
$\mathcal N_{a,E}$ from \eqref{e30}, define  
\[
\Delta_E
:=
P_{\mathcal N_{a,E}^{\perp}}
\mathfrak C\big|_E:
E\longrightarrow\mathcal N_{a,E}^{\perp}.
\]
Every finite-dimensional extension of
$\mathcal N_a$ yields an almost invariant
subspace for the Ces\`aro operator, with
defect bounded above by the dimension
of the added tail space as shown in the following proposition.

\begin{prop}\label{r21}
Let $a\in\mathbb R$. Suppose that $E$ is a finite-dimensional subspace of
$L^2([a,\infty),w(y)\,dy)$. Then
$\mathcal F^{-1}\mathcal N_{a,E}$ is almost invariant for the
Ces\`aro operator $\mathcal C$ on $H^2(\mathbb D)$, and
\begin{equation*}
\operatorname{def}_{\mathcal C}
\bigl(
\mathcal F^{-1}\mathcal N_{a,E}
\bigr)
=
\operatorname{rank}\Delta_E
\leq \dim E.
\end{equation*}
\end{prop}
\begin{proof}
Since $\mathfrak C\mathcal N_a\subseteq\mathcal N_a$,
we have
$P_{\mathcal N_{a,E}^{\perp}}\mathfrak C(h+e)
=\Delta_Ee$
for $h\in\mathcal N_a$ and $e\in E$.
Thus
$P_{\mathcal N_{a,E}^{\perp}}\mathfrak C\mathcal N_{a,E}
=\Delta_EE$, and
$\operatorname{def}_{\mathfrak C}(\mathcal N_{a,E})
=\operatorname{rank}\Delta_E\leq\dim E$.
The result follows by unitary equivalence.
\end{proof}
We continue to use the invariant-subspace lattice described in
\cite[Theorem~6.3]{GallardoPartingtonRoss}   to determine the largest invariant subspace contained in this extension and the smallest invariant subspace that contains it. 
 
\begin{prop}\label{r22}
Let $a\in\mathbb R$ and let $E$ be a finite-dimensional
subspace of $L^2([a,\infty),w(y)\,dy)$. Set
$\mathcal M:=\mathcal F^{-1}\mathcal N_{a,E}$
and $\delta:=e^a$. Then
\(\mathcal I_{\mathcal C}(\mathcal M)
=K_{\Theta_\delta}.
\) In particular, if $E\ne\{0\}$, then
$\operatorname{def}_{\mathcal C}\mathcal M\ge1$.

If $E$ contains a function that is not zero
almost everywhere on any interval $[b,\infty)$ with $b>a$,
then
\(
[\mathcal M]_{\mathcal C}=H^2(\mathbb D).
\)
\end{prop}

\begin{proof}
By \eqref{e17} and the definition of
$\mathcal N_{a,E}$, we have
$K_{\Theta_\delta}\subseteq\mathcal M$ and
$\dim(\mathcal M/K_{\Theta_\delta})=\dim E<\infty$.
Recall that $K_{\Theta_\delta}$ is
$\mathcal C$-invariant. We first prove that
$\mathcal I_{\mathcal C}(\mathcal M)
=K_{\Theta_\delta}$.

By definition,
$\mathcal I_{\mathcal C}(\mathcal M)$
is a closed $\mathcal C$-invariant subspace
satisfying
$K_{\Theta_\delta}\subseteq
\mathcal I_{\mathcal C}(\mathcal M)
\subseteq\mathcal M$.
Consequently,
\[
\dim\bigl(
\mathcal I_{\mathcal C}(\mathcal M)/
K_{\Theta_\delta}
\bigr)
\leq\dim E<\infty.
\]
On the other hand, by
\cite[Theorem~6.3]{GallardoPartingtonRoss},
every proper closed $\mathcal C$-invariant
subspace containing $K_{\Theta_\delta}$
is of the form $K_{\Theta_\beta}$
for some $\beta\ge\delta$.
If $\beta>\delta$, then \eqref{e17}
and the unitarity of $\mathcal F$ in \eqref{e14} give
\[
\dim(K_{\Theta_\beta}/K_{\Theta_\delta})
=\dim L^2((a,\log\beta),w(y)\,dy)=\infty.
\]
Also, $H^2(\mathbb D)/K_{\Theta_\delta}$ is infinite-dimensional by \eqref{e17}.
Thus $K_{\Theta_\delta}$ admits no
proper finite-dimensional closed
$\mathcal C$-invariant extension.
It follows that
$\mathcal I_{\mathcal C}(\mathcal M)
=K_{\Theta_\delta}$.

If $E\ne\{0\}$, then
$K_{\Theta_\delta}\subsetneq\mathcal M$.
Since
$\mathcal I_{\mathcal C}(\mathcal M)
=K_{\Theta_\delta}$,
the space $\mathcal M$ cannot be
$\mathcal C$-invariant. Hence
$\operatorname{def}_{\mathcal C}\mathcal M\ge1$.

Finally, suppose that $E$ contains a function
that is not zero a.e. on any
interval $[b,\infty)$ with $b>a$.
If $[\mathcal M]_{\mathcal C}$ were a proper
subspace of $H^2(\mathbb D)$, then
\cite[Theorem~6.3]{GallardoPartingtonRoss}
would give
$[\mathcal M]_{\mathcal C}=K_{\Theta_\beta}$
for some $\beta\ge\delta$.
Since
$\mathcal M\subseteq[\mathcal M]_{\mathcal C}$,
applying $\mathcal F$ in \eqref{e14} and \eqref{e17} yields
$\mathcal N_a\oplus E
\subseteq\mathcal N_{\log\beta}$.
As every function in $E$ is supported
in $[a,\infty)$, this implies
\[
E\subseteq
L^2((a,\log\beta),w(y)\,dy).
\]
Consequently, every function in $E$
vanishes almost everywhere on
$[b,\infty)$ whenever
$b>\max\{a,\log\beta\}$,
contrary to the hypothesis.
Therefore
$[\mathcal M]_{\mathcal C}=H^2(\mathbb D)$.
\end{proof}

For the special tails $E_{a,\Lambda,\bm m}$, the
corresponding extensions are nearly invariant for
$\{S_t\}_{t\ge0}$ (Theorem~\ref{r13}), but by Propositions~\ref{r14}
and~\ref{r15} they are invariant under neither
$S_t$ nor $S_t^*$ for $t>0$. Applying Proposition~\ref{r21} to
$E_{a,\Lambda,\bm m}$ gives the preliminary defect bound
$\sum_{\nu=1}^r m_\nu$.  We now sharpen this defect. For $\lambda\in\Lambda$, let $m(\lambda)$ denote its
multiplicity, so that $m(\lambda_\nu)=m_\nu$, and set
$m(\lambda)=0$ for $\lambda\notin\Lambda$. Define
\begin{equation*}
E_{a,\Lambda,\bm m}^{+}
:=
\left\{
h\in E_{a,\Lambda,\bm m}:
e^yh\in E_{a,\Lambda,\bm m}
\right\}.
\end{equation*}

Differentiating the integral formula \eqref{e56}
allows us to describe the kernel of the defect map
by studying multiplication by $e^y$. 
\begin{lem}\label{r23} Let $a\in\mathbb R$, 
 $D=\frac{\mathrm d}{\mathrm d y}$ denote ordinary differentiation on $(a,\infty)$,  
$\Lambda=\{\lambda_1,\ldots,\lambda_r\}\subset\mathbb C$
consist of distinct points, and 
$\bm m=(m_1,\ldots,m_r)\in\mathbb N_+^r$. Then
\(
D E_{a,\Lambda,\bm m}
\subseteq
E_{a,\Lambda,\bm m}.
\)
Moreover, the map
\[
\Phi:
\ker\Delta_{E_{a,\Lambda,\bm m}}
\longrightarrow
E_{a,\Lambda,\bm m}^{+},
\qquad
\Phi(g)
=
\chi_{[a,\infty)}\mathfrak Cg,
\]
is a linear isomorphism whose inverse is given by
$\Psi(h) =
\left(\frac12-D\right)h+2e^yh.$  In particular,
\begin{equation}\label{e59}
\dim\ker\Delta_{E_{a,\Lambda,\bm m}}
=
\dim E_{a,\Lambda,\bm m}^{+}.
\end{equation}
\end{lem}

\begin{proof}
For $0\leq k<m_\nu$, differentiation on $(a,\infty)$ gives $Dr_{a,\lambda_\nu,0}
= \lambda_\nu r_{a,\lambda_\nu,0},$ and
$ Dr_{a,\lambda_\nu,k} =\lambda_\nu r_{a,\lambda_\nu,k}+
r_{a,\lambda_\nu,k-1},
\; k\in \mathbb N_+.$ Hence
$DE_{a,\Lambda,\bm m}
\subseteq
E_{a,\Lambda,\bm m}.$ 

Let $g\in\ker\Delta_{E_{a,\Lambda,\bm m}}$ and set
$h:=\chi_{[a,\infty)}\mathfrak Cg$. Since
\[
\Delta_{E_{a,\Lambda,\bm m}}g=0
\quad\Longleftrightarrow\quad
\mathfrak Cg
\in
\mathcal N_a\oplus E_{a,\Lambda,\bm m}\quad\Longleftrightarrow\quad h\in E_{a,\Lambda,\bm m}, \] because 
 $\mathcal N_a$ vanishes on $[a,\infty)$. Formula
\eqref{e56} shows that $h$ is locally
absolutely continuous on $(a,\infty)$.  Differentiating there gives,
almost everywhere,
\begin{equation}\label{e60}
g
=
\left(\frac12-D\right)h+2e^yh.
\end{equation}
Since $g,h,Dh\in E_{a,\Lambda,\bm m}$, equation
\eqref{e60} implies
$e^yh\in E_{a,\Lambda,\bm m}$.  The functions involved are continuous
on $(a,\infty)$, so the almost-everywhere identity has the usual
pointwise representative there. Thus
$ h\in E_{a,\Lambda,\bm m}^{+},$ and hence $\Phi$ is well defined.

Conversely, let $ h\in E_{a,\Lambda,\bm m}^{+}$ and set
$
g := \Psi(h)
= \left(\frac12-D\right)h+2e^yh.
$  Since $ 
Dh\in E_{a,\Lambda,\bm m}$  and  
$e^yh\in E_{a,\Lambda,\bm m},$ one has
$ g\in E_{a,\Lambda,\bm m}. $ Moreover,
\[
D
\left(
e^{-y/2-2e^y}h(y)
\right)
=
-e^{-y/2-2e^y}g(y).
\]
Since $h$ is an exponential polynomial,
$
e^{-y/2-2e^y}h(y)
\longrightarrow0$ as $ y\to\infty.$  Integrating from $y$ to $\infty$ yields
\[
h(y)
=
e^{y/2+2e^y}
\int_y^\infty
e^{-s/2-2e^s}g(s)\,ds
=
(\mathfrak Cg)(y),
\qquad y\geq a.
\]
Therefore $ \mathfrak Cg-h\in\mathcal N_a,$ 
and hence $ \mathfrak Cg
\in \mathcal N_a\oplus E_{a,\Lambda,\bm m}$. Thus
$g\in \ker\Delta_{E_{a,\Lambda,\bm m}}$ 
and  $ \Phi(g)=h.$ Together with
\eqref{e60}, this shows that
$\Phi$ and $\Psi$ are mutually inverse linear maps. Consequently,
\( \dim\ker\Delta_{E_{a,\Lambda,\bm m}} =
\dim E_{a,\Lambda,\bm m}^{+},\)
which proves \eqref{e59}.
\end{proof}

\begin{thm}\label{r24}
Let $a\in\mathbb R$, 
$\Lambda=\{\lambda_1,\ldots,\lambda_r\}\subset\mathbb C$
consist of distinct points, and  
$\bm m=(m_1,\ldots,m_r)\in\mathbb N_+^r$. Then
$\mathcal F^{-1}\mathcal N_{a,\Lambda,\bm m}$ is almost invariant
for the Ces\`aro operator $\mathcal C$, and
\begin{equation}\label{e61}
\operatorname{def}_{\mathcal C}
\bigl(
\mathcal F^{-1}\mathcal N_{a,\Lambda,\bm m}
\bigr)
=
\sum_{\nu=1}^r m_\nu
-
\sum_{\lambda\in\Lambda}
\min\{m(\lambda),m(\lambda+1)\}.
\end{equation}
\end{thm}

\begin{proof}
By Proposition~\ref{r21},
Lemma~\ref{r23}, and the rank--nullity
theorem,
\begin{equation}\label{e62}
\operatorname{def}_{\mathcal C}
\bigl(
\mathcal F^{-1}\mathcal N_{a,\Lambda,\bm m}
\bigr)
=
\sum_{\nu=1}^r m_\nu
-
\dim E_{a,\Lambda,\bm m}^{+}.
\end{equation}
It therefore remains to compute
$\dim E_{a,\Lambda,\bm m}^{+}$.  For $y\geq a$,
\begin{equation}\label{e63}
e^yr_{a,\lambda,k}(y)
=
e^a r_{a,\lambda+1,k}(y),
\end{equation}
which shifts the exponent from $\lambda$
to $\lambda+1$ without changing the polynomial degree.

Given $h\in E_{a,\Lambda,\bm m}$ in \eqref{e47}, we write
$h=\sum_{\nu=1}^r h_\nu$ with $h_\nu\in E_{a,\lambda_\nu,m_\nu}.$ By the linear independence of exponential polynomials with
distinct exponents, the condition
$e^yh\in E_{a,\Lambda,\bm m}$ can be checked blockwise.

If $\lambda_\nu+1\notin\Lambda$, then necessarily $h_\nu=0$.
If $\lambda_\nu+1\in\Lambda$ and  $h_\nu
=
\sum_{k=0}^{m_\nu-1}
c_{\nu,k}r_{a,\lambda_\nu,k},$ 
then \eqref{e63} gives
\[
e^yh_\nu
=
e^a
\sum_{k=0}^{m_\nu-1}
c_{\nu,k}r_{a,\lambda_\nu+1,k}.
\]
Hence,
\(
e^yh_\nu
\in
E_{a,\lambda_\nu+1,m(\lambda_\nu+1)}\) if and only if $
c_{\nu,k}=0,\;
k\geq m(\lambda_\nu+1).$  Therefore, the part of $E_{a,\lambda_\nu,m_\nu}$ contributing to
$E_{a,\Lambda,\bm m}^{+}$ is precisely
\[
\bigvee
\left\{
r_{a,\lambda_\nu,k}:
0\leq k<
\min\{m_\nu,m(\lambda_\nu+1)\}
\right\},
\]
and thus has dimension
$ \min\{m_\nu,m(\lambda_\nu+1)\}$ by Lemma \ref{r12}. Consequently,
\begin{equation}\label{e64}
\dim E_{a,\Lambda,\bm m}^{+}
=
\sum_{\lambda\in\Lambda}
\min\{m(\lambda),m(\lambda+1)\}.
\end{equation}
Combining
\eqref{e62}
and \eqref{e64} gives
\eqref{e61}.
\end{proof}

\begin{rem}\label{r25}
Equivalently, the formula in Theorem~\ref{r24}
can be written as
\begin{equation*} 
\operatorname{def}_{\mathcal C}
\bigl(
\mathcal F^{-1}\mathcal N_{a,\Lambda,\bm m}
\bigr)
=
\sum_{\lambda\in\Lambda}
\max\{m(\lambda)-m(\lambda+1),0\}.
\end{equation*}
Thus the Ces\`aro defect is independent of the support parameter
$a\in\mathbb R$ and depends only on the exponent configuration and its
multiplicities. Since these tails are nonzero exponential polynomials
on every terminal half-line, Proposition~\ref{r22} also shows that
$K_{\Theta_{e^a}}$ is the \emph{largest contained closed
$\mathcal C$-invariant subspace}, while the smallest closed
$\mathcal C$-invariant subspace containing the extension is $H^2$.
\end{rem}

As an application of the general defect formula,
we now consider the smallest nearly invariant
subspace generated by a single exponential--polynomial
tail. The main point is to establish its minimality
through backward translations.

\begin{thm}\label{r20}
Let $a\in\mathbb R$, $\lambda\in\mathbb C$ and $m\in \mathbb{N}_+$.
The smallest closed subspace containing
$r_{a,\lambda,m-1}$ in \eqref{e34} and nearly invariant
for $\{S_t\}_{t\ge0}$ is
$\mathcal N_a\oplus E_{a,\lambda,m}$,
where $E_{a,\lambda,m}$ is defined in \eqref{ealm}.
The corresponding subspace
$M_{a,\lambda,m}:=
\mathcal F^{-1}(\mathcal N_a\oplus E_{a,\lambda,m})$
of $H^2(\mathbb D)$ satisfies
$\operatorname{def}_{\mathcal C}M_{a,\lambda,m}=m$.
\end{thm}

\begin{proof}
By Theorem~\ref{r13},
$\mathcal N_a\oplus E_{a,\lambda,m}$
is nearly invariant for $\{S_t\}_{t\ge0}$
and contains $r_{a,\lambda,m-1}$.
It remains to prove that this is the
smallest such closed subspace.

Let $M$ be any closed nearly invariant
subspace containing $r_{a,\lambda,m-1}$.
By Lemma~\ref{r11}, for every $t\ge0$,
the function $r_{a,\lambda,m-1}(\cdot+t)$
belongs to $L_w^2$. Moreover,
\[
S_t\bigl(r_{a,\lambda,m-1}(\cdot+t)\bigr)
=r_{a,\lambda,m-1}\in M.
\]
Since $M$ is nearly invariant, we obtain
\(r_{a,\lambda,m-1}(\cdot+t)\in M\) for $t\ge0.$

We show that
$\mathcal N_a\oplus E_{a,\lambda,m}\subseteq M$.
Let $h\in M^\perp$. Then
$\langle r_{a,\lambda,m-1}(\cdot+t),h\rangle_{L_w^2}=0$
for $t\ge0$. By the definition of $r_{a,\lambda,m-1}$, it holds that
\begin{align}
0= e^{-\lambda t}
\langle r_{a,\lambda,m-1}(\cdot+t),h
\rangle_{L_w^2}=
\frac{1}{(m-1)!}
\int_{a-t}^{\infty}
(y+t-a)^{m-1}e^{\lambda(y-a)}
\overline{h(y)}w(y)\,dy.\label{0dr}
\end{align}
Differentiating the integral with respect
to $t$, we obtain, for $0\le j\le m-1$,
\[
0= \frac{{\rm d}^j}{{\rm d}t^j}
\left(
e^{-\lambda t}
\langle r_{a,\lambda,m-1}(\cdot+t),h
\rangle_{L_w^2}
\right)=
\frac{1}{(m-1-j)!}
\int_{a-t}^{\infty}
(y+t-a)^{m-1-j}
e^{\lambda(y-a)}
\overline{h(y)}w(y)\,dy.
\]
Taking $j=m-1$, the integral on the right has no polynomial factor. One more differentiation, using the fundamental theorem of calculus for locally integrable functions, yields
\[
0= \frac{{\rm d}^m}{{\rm d}t^m}
\left(
e^{-\lambda t}
\langle r_{a,\lambda,m-1}(\cdot+t),h
\rangle_{L_w^2}
\right)=
e^{-\lambda t}
\overline{h(a-t)}w(a-t),
\qquad \text{for a.e. } t> 0.
\]
 Since $w(y)>0$ for every $y\in\mathbb R$, 
we conclude that
$h(y)=0$ for a.e. $y<a.$ Thus $h\perp\mathcal N_a$.

We next show that $h\perp E_{a,\lambda,m}$.
Since $h=0$ almost everywhere on $(-\infty,a)$,
equation \eqref{0dr} gives
\[
0=\frac{1}{(m-1)!}
\int_a^\infty (y+t-a)^{m-1}
e^{\lambda(y-a)}\overline{h(y)}w(y)\,dy=\sum_{k=0}^{m-1}
\frac{t^{m-1-k}}{(m-1-k)!}
\langle r_{a,\lambda,k},h\rangle_{L_w^2},
\qquad t\ge0,
\]
where the second equality follows from the
binomial formula and the definition of
$r_{a,\lambda,k}$.  The right-hand side is a polynomial in $t$
that vanishes for every $t\ge0$.
Hence all its coefficients are zero, that is,
$\langle r_{a,\lambda,k},h\rangle_{L_w^2}=0$
for $0\le k<m$. Since the definition of
$E_{a,\lambda,m}$ in \eqref{ealm},
we obtain $h\perp E_{a,\lambda,m}$.

Thus every $h\in M^\perp$ is orthogonal to
both $\mathcal N_a$ and $E_{a,\lambda,m}$.
Consequently,
\[
M^\perp\subseteq
(\mathcal N_a\oplus E_{a,\lambda,m})^\perp,
\]
and taking orthogonal complements yields
$\mathcal N_a\oplus E_{a,\lambda,m}\subseteq M$.
Since $\mathcal N_a\oplus E_{a,\lambda,m}$ is closed, nearly invariant for $\{S_t\}_{t\ge0}$, and contains $r_{a,\lambda,m-1}$, the claimed minimality follows.

Finally, the Ces\`aro defect follows
from Theorem~\ref{r24}. Taking
$\Lambda=\{\lambda\}$ with multiplicity $m$,
we have $m(\lambda+1)=0$, because
$\lambda+1\notin\Lambda$. Consequently,
$\operatorname{def}_{\mathcal C}M_{a,\lambda,m}
=m-\min\{m,m(\lambda+1)\}=m.$\end{proof}

\begin{cor}
For the space $M_{a,\lambda,m}$ in
Theorem~\ref{r20},
$ 
\{f\in M_{a,\lambda,m}:\mathcal Cf\in M_{a,\lambda,m}\}
=K_{\Theta_{e^a}},$ i.e., $K_{\Theta_{e^a}}$ is the largest
closed $\mathcal C$-invariant subspace
contained in $M_{a,\lambda,m}$.
\end{cor}

\begin{proof}
By Proposition~\ref{r21} and Theorem~\ref{r20},
$\operatorname{rank}\Delta_{E_{a,\lambda,m}}
=m=\dim E_{a,\lambda,m}$;
hence $\Delta_{E_{a,\lambda,m}}$ is injective.

Suppose that $f\in M_{a,\lambda,m}$ and
$\mathcal Cf\in M_{a,\lambda,m}$.
Write $\mathcal Ff=h+e$, where
$h\in\mathcal N_a$ and
$e\in E_{a,\lambda,m}$.
Since $\mathfrak C\mathcal N_a\subseteq\mathcal N_a$,
we have
\[
\mathfrak Ce
=\mathcal F\mathcal Cf-\mathfrak Ch
\in\mathcal N_a\oplus E_{a,\lambda,m}.
\]
Thus $\Delta_{E_{a,\lambda,m}}e=0$.
By injectivity, $e=0$, and hence
$f\in\mathcal F^{-1}\mathcal N_a
=K_{\Theta_{e^a}}$.

The reverse inclusion follows from the
$\mathcal C$-invariance of $K_{\Theta_{e^a}}$.
Moreover, every $\mathcal C$-invariant
subspace contained in $M_{a,\lambda,m}$
must be contained in $K_{\Theta_{e^a}}$,
which proves maximality.
\end{proof}

For the model-space family in Theorem~\ref{r17},
we can further determine the exact Ces\`aro
defect and its corresponding defect space.

\begin{cor}\label{r26}
Let $N\in \mathbb{N}_+$ and assume the parameters in
Corollary~\ref{r19}. Then
\begin{equation}\label{e65}
\operatorname{def}_{\mathcal C}
(K_{z^N\Theta_\delta})
=
1.
\end{equation}

Moreover, the one-dimensional Ces\`aro defect space is
\begin{equation}\label{e66}
P_{K_{z^N\Theta_\delta}^{\perp}}
\mathcal C K_{z^N\Theta_\delta}
= \mathbb C\,
P_{z^N\Theta_\delta H^2}
(\mathcal C\Theta_\delta).
\end{equation}
\end{cor}

\begin{proof}
The set $\Lambda_N$ is a single consecutive chain of length $N$
with all multiplicities equal to one. Hence
Theorem~\ref{r24} gives
\eqref{e65}.  To identify the  corresponding defect space, set
\[
q_{a,j}(y)
:=
e^{(j+\frac12)y}\chi_{[a,\infty)}(y),
\qquad
0\leq j\leq N-1.
\]
Since each $q_{a,j}$ is a nonzero scalar multiple of
$r_{a,j+\frac12,0}$,
\[
E_{a,\Lambda_N,\bm1}
=
\bigvee
\left\{
q_{a,j}:0\leq j\leq N-1
\right\}.
\]
Moreover,
\[
E_{a,\Lambda_N,\bm1}^{+}
=
\bigvee
\left\{
q_{a,j}:0\leq j\leq N-2
\right\},
\]
where the latter span is understood to be $\{0\}$ when $N=1$.

For $0\leq j\leq N-2$, the inverse map $\Psi$ in
Lemma~\ref{r23} satisfies
\[
\Psi(q_{a,j})
=
\left(\frac12-D\right)q_{a,j}
+2e^yq_{a,j}
=
-jq_{a,j}+2q_{a,j+1}.
\]
These vectors form a triangular family with respect to
$q_{a,1},\ldots,q_{a,N-1}$. Hence
\[
\ker\Delta_{E_{a,\Lambda_N,\bm1}}
=
\bigvee
\left\{
q_{a,j}:1\leq j\leq N-1
\right\},
\]
with the right-hand side understood to be $\{0\}$ when $N=1$.
In particular, $q_{a,0}
\notin
\ker\Delta_{E_{a,\Lambda_N,\bm1}}.$ Since
$\operatorname{rank}\Delta_{E_{a,\Lambda_N,\bm 1}}=1$, it follows
that  $\operatorname{ran}
\Delta_{E_{a,\Lambda_N,\bm1}}
= \mathbb C\,
\Delta_{E_{a,\Lambda_N,\bm1}}q_{a,0}.$

A direct calculation from \eqref{e11}--\eqref{e14} gives
$\mathcal F^{-1}q_{a,0}=
(e^{1-\delta}/\sqrt2)\Theta_\delta$.
Since $K_{z^N\Theta_\delta}^{\perp}
=z^N\Theta_\delta H^2$, unitary equivalence yields
\[
P_{K_{z^N\Theta_\delta}^{\perp}}
\mathcal C K_{z^N\Theta_\delta}
=
\mathbb C\,
P_{z^N\Theta_\delta H^2}
(\mathcal C\Theta_\delta),
\]
which proves \eqref{e66}.
\end{proof}

In summary, the Jordan structure governs the backward-translation
dynamics, while the successor relation $\lambda\mapsto\lambda+1$
governs the Ces\`aro defect. For the distinguished model-space
family, the exponents form a single consecutive chain, so the
Ces\`aro defect is one for every $N\in \mathbb{N}_+$.
\subsection{Finite Blaschke extensions and semigroup near invariance}
We now turn to our main questions on the finite Blaschke extensions of $K_{\Theta_\delta}$.
The decomposition \eqref{e24} reduces the calculation of their
Ces\`aro defects to the finite-dimensional space $E_{\delta,B}$.
We determine the defect of $\mathcal C$ and then characterize which extensions are nearly invariant for the composition semigroup.

\begin{thm}\label{r27}
Let $\delta>0$ and $d\in\mathbb N_+$. Let $B$ be the finite
Blaschke product in \eqref{e22}, whose half-plane counterpart
$B^+$ has distinct zeros $\alpha_1,\ldots,\alpha_d$ of
multiplicities $n_1,\ldots,n_d$. Then
$\operatorname{def}_{\mathcal C}K_{B\Theta_\delta}=d$.
\end{thm}

\begin{proof}
Put $N:=\sum_{j=1}^d n_j$. By \eqref{e24},
\[
\mathcal J K_{B\Theta_\delta}
=L^2(0,\delta)\oplus E_{\delta,B},
\qquad \dim E_{\delta,B}=N.
\]
Recall that $V=\mathcal J\mathcal C\mathcal J^{-1}$ and that
$L^2(0,\delta)$ is invariant under $V$. Define
\[
\Gamma_{\delta,B}
:=P_{(\mathcal J K_{B\Theta_\delta})^\perp}
  V\big|_{E_{\delta,B}}.
\]
Since $\mathcal J$ is unitary and $L^2(0,\delta)$ is
$V$-invariant, the defect is determined by the action
of $V$ on the $N$-dimensional space $E_{\delta,B}$. The rank--nullity theorem gives
\begin{equation}\label{eq:r27-rank}
\operatorname{def}_{\mathcal C}K_{B\Theta_\delta}
=\operatorname{rank}\Gamma_{\delta,B}
=N-\dim\ker\Gamma_{\delta,B}.
\end{equation}
We therefore compute $\ker\Gamma_{\delta,B}$. For $f\in E_{\delta,B}$,
\begin{equation}\label{equivs}
f\in\ker\Gamma_{\delta,B}
 \iff Vf\in L^2(0,\delta)\oplus E_{\delta,B} 
 \iff \chi_{[\delta,\infty)}Vf\in E_{\delta,B}.
\end{equation}

First suppose that $f\in\ker\Gamma_{\delta,B}$ and set
$h:=\chi_{[\delta,\infty)}Vf$. Then $h\in E_{\delta,B}$, so
\begin{equation}\label{hfp}
h(x)=\sum_{j=1}^d p_j(x-\delta)
 e^{-\overline{\alpha_j}(x-\delta)},\qquad x>\delta,
\end{equation}
where $\deg p_j\le n_j-1$. Put
$\mathscr D:=x(1-\frac{\mathrm d}{\mathrm dx})$.
Since $h=Vf$ on $(\delta,\infty)$, the differential identity
\eqref{e21} gives $f=\mathscr D h$ a.e. on this interval.

Writing $u=x-\delta$, we obtain
\begin{equation}\label{eq:r27-differential}
\mathscr D\bigl(p_j(u)e^{-\overline{\alpha_j}u}\bigr)
=(\delta+u)\bigl((1+\overline{\alpha_j})p_j(u)-p_j'(u)\bigr)
 e^{-\overline{\alpha_j}u},
\end{equation}
where  $u=x-\delta$.
If $p_j\ne0$ has degree $k$ and leading coefficient $c_k\ne0$,
the polynomial factor on the right has leading term
$(1+\overline{\alpha_j})c_k u^{k+1}$. This coefficient is nonzero
because $\operatorname{Re}\alpha_j>0$. Thus applying $\mathscr D$
raises the polynomial degree in each nonzero exponential block by
exactly one. Exponential polynomials with distinct rates are
linearly independent, so $f=\mathscr D h\in E_{\delta,B}$ forces
$\deg p_j+1\le n_j-1$ for every nonzero $p_j$. Consequently,
$\deg p_j\le n_j-2$ for $1\le j\le d$, with $p_j=0$ when $n_j=1$.
This proves the necessary condition.

Conversely, choose any polynomials $p_j$ with
$\deg p_j\le n_j-2$, taking $p_j=0$ if $n_j=1$.
Define $h$ on $(\delta,\infty)$ by \eqref{hfp} and put
$f:=\chi_{[\delta,\infty)}\mathscr D h$, where differentiation is
performed only on $(\delta,\infty)$. Formula
\eqref{eq:r27-differential} shows that $f\in E_{\delta,B}$.
It remains to check that $f$ belongs to the kernel.
For $x>\delta$, we have $f(x)=x(h(x)-h'(x))$ and hence
$e^{-x}f(x)/x=-\frac{\mathrm d}{\mathrm dx}(e^{-x}h(x))$.
Since $h$ is a finite sum of polynomial--exponential terms with
$\operatorname{Re}\alpha_j>0$, we have
$\lim_{x\to\infty}e^{-x}h(x)=0$.
Using the integral formula \eqref{e20}, we obtain, for $x>\delta$,
\[
(Vf)(x)
=e^x\int_x^\infty\frac{e^{-s}}{s}f(s)\,\mathrm ds
=-e^x\int_x^\infty
  \frac{\mathrm d}{\mathrm ds}\bigl(e^{-s}h(s)\bigr)\,\mathrm ds
=h(x).
\]
Thus $\chi_{[\delta,\infty)}Vf=h\in E_{\delta,B}$, and
\eqref{equivs} implies $f\in\ker\Gamma_{\delta,B}$.
We have proved both directions and therefore obtained
\begin{equation}\label{eq:r27-kernel}
\ker\Gamma_{\delta,B}
=\{\chi_{[\delta,\infty)}(x)
\sum_{j=1}^d\mathscr D\bigl(
 p_j(x-\delta)e^{-\overline{\alpha_j}(x-\delta)}\bigr):  \deg p_j\le n_j-2,\quad 1\le j\le d\}.
\end{equation}
As above, $p_j=0$ if $n_j=1$, and every derivative in
\eqref{eq:r27-kernel} is taken on $(\delta,\infty)$.

Finally, for each $1\le j\le d$, consider the linear map
\[
A_j:\mathcal P_{n_j-2}\longrightarrow\mathcal P_{n_j-1},
\qquad
A_jp(u):=(\delta+u)
\bigl((1+\overline{\alpha_j})p(u)-p'(u)\bigr),
\]
where $\mathcal P_m$ denotes the space of polynomials
of degree at most $m$, with $\mathcal P_{-1}=\{0\}$.
By \eqref{eq:r27-differential}, every nonzero
$p\in\mathcal P_{n_j-2}$ satisfies
$\deg(A_jp)=\deg p+1$.
Thus $A_j$ is injective and
$\dim\operatorname{ran}A_j=n_j-1$.
Since exponential polynomials with distinct rates
are linearly independent, \eqref{eq:r27-kernel} yields
\[
\dim\ker\Gamma_{\delta,B}
=\sum_{j=1}^d\dim\operatorname{ran}A_j
=\sum_{j=1}^d(n_j-1)=N-d.
\]
The conclusion \(\operatorname{def}_{\mathcal C}K_{B\Theta_\delta}
=N-(N-d)=d\) now follows from \eqref{eq:r27-rank}.
\end{proof}

\begin{exm}\label{r29}
Write $\mathfrak b_a(z):=(z-a)/(1-\overline a z)$.
For   $N\in \mathbb{N}_+$ and $B(z)=z^N$, all half-plane zeros equal $1$, so
\(\operatorname{def}_{\mathcal C}K_{z^N\Theta_\delta}=1.
\) 
For $B(z)=z^2\mathfrak b_{1/2}(z)$, the corresponding
half-plane zeros are $\alpha_1=1$ and $\alpha_2=3$,
with multiplicities $2$ and $1$, respectively.
 Hence
\(
\operatorname{def}_{\mathcal C}K_{B\Theta_\delta}=2.
\)
These examples show that the Ces\`aro defect
depends on the number of distinct half-plane
zeros rather than their multiplicities.
\end{exm}

We next characterize the finite Blaschke extensions
of $K_{\Theta_\delta}$ that are also nearly invariant
for the associated composition semigroup
$\{C_t\}_{t\geq0}$.

\begin{thm}\label{r30}
Let $\delta>0$ and let $B$ be a nonconstant finite
Blaschke product. Then $K_{B\Theta_\delta}$ is nearly
invariant for the composition semigroup
$\{C_t\}_{t\ge0}$ if and only if
$B(z)=z^N$ for some $N\ge1$, up to a
unimodular constant.
\end{thm}

\begin{proof}
\emph{ Sufficiency}. Suppose that
$B(z)=z^N$ up to a unimodular constant,
where $N\ge1$. Since a unimodular constant
does not change the corresponding model
space, we may assume that $B(z)=z^N$.

Put $a=\log\delta$ and
$\Lambda_N=\{\frac12,\frac32,\ldots,N-\frac12\}$.
By Corollary~\ref{r19},
\[
\mathcal F K_{z^N\Theta_\delta}
=\mathcal N_a\oplus E_{a,\Lambda_N,\bm1}.
\]
Theorem~\ref{r13} shows that this space
is nearly invariant for $\{S_t\}_{t\ge0}$.
Since
$\mathcal F C_t\mathcal F^{-1}=e^{t/2}S_t$,
the space $K_{z^N\Theta_\delta}$ is nearly
invariant for $\{C_t\}_{t\ge0}$.

\emph{Necessity}. Suppose that $K_{B\Theta_\delta}$
is nearly invariant for $\{C_t\}_{t\ge0}$.
Let $\alpha_1,\ldots,\alpha_d\in\mathbb C_+$
be the distinct half-plane zeros of $B^+$,
with multiplicities $n_1,\ldots,n_d$.
Put $a:=\log\delta$ and $c_j:=1-\overline{\alpha_j}$, so
$\operatorname{Re}c_j=1-\operatorname{Re}\alpha_j<1$ for each $j$.
Recall that $\mathcal F=WT^{-1}\mathcal J$ in \eqref{e14}.
By the definitions of $W$ and $T$, we have
\[
WT^{-1}\left(
(x-\delta)^k e^{-\overline{\alpha_j}(x-\delta)}
\chi_{[\delta,\infty)}(x)\right)(y)=e^{\overline{\alpha_j}\delta-1}
e^{y/2}e^{c_j e^y}(e^y-\delta)^k
\chi_{[a,\infty)}(y).
\]
Then  by \eqref{e24},
$\mathcal F K_{B\Theta_\delta}
=\mathcal N_a\oplus E,$ where
\[
E:=\bigvee\left\{v_{j,k}:1\le j\le d,\ 0\le k<n_j\right\},\qquad
v_{j,k}(y):=
e^{y/2}e^{c_j e^y}(e^y-\delta)^k
\chi_{[a,\infty)}(y).
\]
Since $\mathcal F C_t\mathcal F^{-1}
=e^{t/2}S_t$, the space
$\mathcal N_a\oplus E$ is nearly invariant
for $\{S_t\}_{t\ge0}$.

To prove that $c_j=0$ for every $j$,
it suffices to consider the functions
$v_{j,0}$, which contain all the distinct
exponential parameters $c_j$. Fix $1\le j\le d$. By the domain
characterization of $\mathcal B_t$,
\[
\begin{aligned}
v_{j,0}\in D(\mathcal B_t)
&\Longleftrightarrow
\int_{\mathbb R}|v_{j,0}(y)|^2
w(y-t)\,dy<\infty\\
&\Longleftrightarrow
\int_\delta^\infty
e^{2(\operatorname{Re}c_j-e^{-t})x}\,dx
<\infty\\
&\Longleftrightarrow
\operatorname{Re}c_j<e^{-t}.
\end{aligned}
\]
Since $\operatorname{Re}c_j<1$, there exists
$\varepsilon>0$ such that
$\operatorname{Re}c_j<e^{-t}$
for every $1\le j\le d$ and
$0<t<\varepsilon$.
Hence $v_{j,0}\in D(\mathcal B_t)$,
or equivalently,
$v_{j,0}(\cdot+t)\in L_w^2(\mathbb R)$,
for all $0<t<\varepsilon$.   Now
$S_t(v_{j,0}(\cdot+t))=v_{j,0}\in E$.
By near invariance,
$v_{j,0}(\cdot+t)\in\mathcal N_a\oplus E$.
Taking its restriction to $[a,\infty)$,
we obtain
\[
e^{t/2}e^{y/2}e^{c_j e^t e^y}
\chi_{[a,\infty)}(y)\in E,
\qquad 0<t<\varepsilon.
\]
By the definition of $E$, there exist
polynomials $p_{\ell,t}$ of degree
less than $n_\ell$ such that
\begin{align}
e^{c_j e^t x}
=\sum_{\ell=1}^d
p_{\ell,t}(x-\delta)e^{c_\ell x},
\qquad x>\delta,\label{es}
\end{align}
where the nonzero factor $e^{t/2}$
has been absorbed into the coefficients.
The equality holds everywhere on
$(\delta,\infty)$ by continuity.

We claim that
$c_j e^t\in\{c_1,\ldots,c_d\}$
for every $0<t<\varepsilon$.
Indeed, applying the differential operator
$\prod_{\ell=1}^d(D-c_\ell)^{n_\ell}$,
where $D=\frac{\rm d}{{\rm d}x}$, to \eqref{es}
gives
\[
\left(
\prod_{\ell=1}^d
(c_j e^t-c_\ell)^{n_\ell}
\right)e^{c_j e^t x}=0.
\]
The right-hand side of the preceding
exponential-polynomial identity is
annihilated by this operator because
$\deg p_{\ell,t}<n_\ell$.
Since $e^{c_j e^t x}\ne0$, we obtain
$\prod_{\ell=1}^d(c_j e^t-c_\ell)^{n_\ell}=0$,
and hence
$c_j e^t\in\{c_1,\ldots,c_d\}$ for $0<t<\varepsilon$.
Therefore $c_j e^t$ must be constant
on the interval $(0, \varepsilon)$.
Letting $t\downarrow0$, we obtain
$c_j e^t=c_j$ for $0<t<\varepsilon$.
Since $e^t\ne1$ for $t>0$, it follows
that $c_j=0$. As $j$ was arbitrary, all $c_j$ vanish.
Thus $\alpha_j=1$ for every $j$.
Under the Cayley transform, the point
$\alpha_j=1$ corresponds to the disk
point $z=0$. Hence all zeros of $B$
are at the origin, and therefore
$B(z)=z^N$ for some $N\ge1$,
up to a unimodular constant.
\end{proof}
 
We conclude this section by summarizing how the preceding results answer the main questions raised in the Introduction.
\begin{rem}
  For a nonconstant finite Blaschke product
$B$ with $d$ distinct zeros, Theorem~\ref{r27} and
Proposition~\ref{r22} give
\[
\operatorname{def}_{\mathcal C}K_{B\Theta_\delta}=d,
\qquad
\mathcal I_{\mathcal C}(K_{B\Theta_\delta})
=K_{\Theta_\delta},
\qquad
[K_{B\Theta_\delta}]_{\mathcal C}=H^2.
\]
Thus the defect depends only on the number of distinct
zeros of $B$, not on their multiplicities. Moreover,
$K_{B\Theta_\delta}$ is almost, but not strictly,
$\mathcal C$-invariant, with $K_{\Theta_\delta}$ as its
largest contained and $H^2$ as its smallest containing
closed $\mathcal C$-invariant subspace.
\end{rem}

\section*{Declarations}

\noindent\textbf{Acknowledgments.}
This work was supported by the National Natural Science Foundation of China (Grant No. 12471126).  The author would like to thank Professor Jonathan R. Partington (University of Leeds) for his helpful discussions and valuable suggestions during the preparation of this paper.

\medskip
\noindent\textbf{Data availability.}
Not applicable.

\medskip
\noindent\textbf{Conflict of interest.}
The author declares no relevant financial or non-financial interests.

\bigskip
\begingroup\small
  \noindent\textbf{ Yuxia Liang }\par \vspace{1mm}
 \noindent {School of Mathematical Sciences, Tianjin Normal University, 
Tianjin 300387, P.R. China }\par \vspace{1mm}
\noindent\textit{E-mail address:} {liangyx1986@126.com} \par
\endgroup
\end{document}